\documentclass[amstex,12pt, amssymb]{article}

\usepackage{mathtext}
\usepackage[cp1251]{inputenc}
\usepackage[T2A]{fontenc}
\usepackage[dvips]{graphicx}
\usepackage{amsmath}
\usepackage{amssymb}
\usepackage{amsxtra}
\usepackage{latexsym}
\usepackage{ifthen}

\newcounter{lemma}[section]

\newcounter{corollary}[section]

\newcounter{remark}[section]

\newcounter{theorem}[section]

\newcounter{proposition}[section]

\newcounter{example}

\numberwithin{equation}{section}

\begin{document}

\markboth{V. DESYATKA, E.~SEVOST'YANOV}{\centerline{ON BOUNDARY
H\"{O}LDER CONTINUITY...}}

\def\cc{\setcounter{equation}{0}
\setcounter{figure}{0}\setcounter{table}{0}}

\overfullrule=0pt

%\normalsize\large

\author{VICTORIA DESYATKA, ZARINA KOVBA, EVGENY SEVOST'YANOV}

\title{
{\bf ON BOUNDARY H\"{O}LDER CONTINUITY OF SOBOLEV AND ORLICZ-SOBOLEV
CLASSES}}

\date{\today}
\maketitle

%\large
\begin{abstract}
We investigate distortion estimates for mappings at boundary points
of a domain. We consider mappings of the Sobolev and Orlicz-Sobolev
classes and some other classes of mappings that do not preserve the
boundary of a domain. For above mappings, we establish distortion
estimates at boundary points. In particular, under certain
conditions on the characteristics of the mappings, we show that they
are H\"{o}lder continuous. In the manuscript we consider both the
case of ``good'' boundaries and ``domains with prime ends''. We have
obtained not only H\"{o}lder-type estimates, but also some more
general ones under appropriate (more general) conditions on the
characteristic.
\end{abstract}


\section{Introduction}

This paper is devoted to the study of mappings with bounded and
finite distortion, see, e.g., \cite{A}, \cite{Cr}, \cite{MRV$_1$},
\cite{MRSY}, \cite{PSS}, \cite{RV}, \cite{SalSt}, \cite{Vu}
and~\cite{Va}. In~\cite{DS} we established the possibility of a
continuous boundary extension of open discrete mappings satisfying
Poletskii inequality, which, generally speaking, do not preserve the
boundary of a domain. The present paper continues the research in
this direction. In particular, we establish here that the distortion
of the distance under similar mappings at boundary points can be
obtained in a more or less explicit form.

\medskip
Research in this direction is closely related to some
classical mapping classes. Recall that, any quasiconformal
(quasiregular) mapping is H\"{o}lder continuous inside a domain with
some exponent. In particular, the following result holds (see, e.g.,
\cite[Theorem~3.2]{MRV$_1$}).

\medskip
{\bf Theorem~A. (Martio, Rickman, V\"{a}is\"{a}l\"{a}, 1970)}.
{\it\, Suppose that $f$ is bounded and quasiregular in a domain
$G\subset {\Bbb R}^n$ and $F$ is a compact subset of $G.$ Let
$K_I(f)$ be the smallest constant $K$ for which the inequality $J(x,
f)\leqslant K\cdot (l(f^{\,\prime}(x)))^n$ holds for almost all
$x\in G,$ $J(x, f)=\det f^{\,\prime}(x)$ and
$l(f^{\,\prime}(x))=\min\limits_{|h|=1}|f^{\,\prime}(x)h|.$ Then
there is some constant $\lambda_n$ depending only on $n$ such that
the relation
$$|f(x)-f(y)|\leqslant C|x-y|^{\alpha}$$
holds for $x\in F,$ $y\in G,$ where
$\alpha=(K_I(f))^{\frac{1}{1-n}}$ and $C=\lambda_n(d(f, \partial
G))^{\,-\alpha}d(fG).$}

\medskip
Our immediate goal is to obtain some analogue of this result for
some more general class of mappings at the boundary points of the
domain. The corresponding statement is formulated both for ``good''
boundaries and for domains assuming the use of prime ends. Similar
results have been recently obtained by the third co-author (see,
e.g., \cite{RSS}, \cite{MSS}, \cite{Sev$_4$}, \cite{DovSev} and
\cite{SBDI}). In these papers we somehow assumed that the mapping
preserves the boundary of the domain.

\medskip
Recall some definitions. Let $U$ be an open set in ${\Bbb R}^n.$ In
what follows, $C^k_0(U)$ denotes the space of functions $u:U
\rightarrow {\Bbb R} $ with a compact support in $U,$ having $k$
partial derivatives with respect to any variable that are continuous
in $U.$ Let $U$ be an open set, $U\subset{\Bbb R}^n,$ and let
$u:U\rightarrow {\Bbb R}$ be some function, $u \in L_{\rm
loc}^{\,1}(U).$ Assume that, there exists a function $v\in L_{\rm
loc}^{\,1}(U)$ such that
$$\int\limits_U \frac{\partial \varphi}{\partial x_i}(x)u(x)\,dm(x)=
-\int\limits_U \varphi(x)v(x)\,dm(x)$$
for any function $\varphi\in C_1^{\,0}(U).$ Then we say that the
function $v$ is a weak derivative of the first order of~$u$ with
respect to $x_i$ and denoted by the $\frac{\partial u}{\partial
x_i}(x):= v.$ We write $u\in W_{\rm loc}^{1,1}(U)$ if $u$ has weak
derivatives of the first order with respect to each of the variables
in $U$ and $\frac{\partial u}{\partial x_i}$ are locally integrable
in~$U.$

\medskip
A mapping $f:U\rightarrow {\Bbb R}^n$ belongs to the Sobolev class
$W_{\rm loc}^{1,1}(U),$ we write $f \in W^{1,1}_{\rm loc}(U),$ if
all coordinate functions of $f=(f_1,\ldots, f_n)$ have weak partial
derivatives of the first order, which are locally integrable in $U.$
We write $f\in W^{1, p}_{\rm loc}(U),$ $p\geqslant 1,$ if all
coordinate functions $f_i,$ $1\leqslant i\leqslant n,$ of
$f=(f_1,\ldots, f_n)$ have weak partial derivatives of the first
order, which are locally integrable in $U$ to the degree $p.$

\medskip
Let $\varphi:{\Bbb R}^+\rightarrow{\Bbb R}^+$ be a measurable
function. The {\it Orlicz-Sobolev class} $W^{1,\varphi}_{\rm
loc}(D)$ is the class of all locally integrable mappings $f$ with
the first distributional derivatives whose gradient $\nabla f$
belongs locally in $D$ to the Orlicz class. By the definition,
$W^{1,\varphi}_{\rm loc}\subset W^{1,1}_{\rm loc}$. For the case
when $\varphi(t)=t^p$, $p\geqslant 1,$ we write as usual $f\in
W^{1,p}_{\rm loc}.$

\medskip
In more detail, we write $f\in W^{1,\varphi}_{\rm loc}(D),$ if
$f_i\in W^{1,1}_{\rm loc}(D)$ and
$$
\int\limits_{G}\varphi\left(|\nabla f(x)|\right)\,dm(x)<\infty
$$
for any domain $G\subset D$ with $\overline{G}\subset D,$ $|\nabla
f(x)|=\sqrt{\sum\limits_{i,j}\left(\frac{\partial f_i}{\partial
x_j}\right)^2}.$

\medskip
Assume that, a mapping $f$ has partial derivatives almost everywhere
in $D.$ In this case, we set
\begin{gather*}l\left(f^{\,\prime}(x)\right)\,=\,\min\limits_{h\in {\Bbb R}^n
\backslash \{0\}} \frac {|f^{\,\prime}(x)h|}{|h|}\,, \label{eq5_a}
\,\Vert f^{\,\prime}(x)\Vert\,=\,\max\limits_{h\in {\Bbb
R}^n \backslash \{0\}} \frac {|f^{\,\prime}(x)h|}{|h|}\,,\\
J(x,f)=\det f^{\,\prime}(x)\,,\nonumber\end{gather*}
and define for any $x\in D$ and $\alpha\geqslant 1$
\begin{equation}\label{eq0.1.1A}
K_I(x,f)\quad =\quad\left\{
\begin{array}{rr}
\frac{|J(x,f)|}{{l\left(f^{\,\prime}(x)\right)}^{\,n}}, & J(x,f)\ne 0,\\
1,  &  f^{\,\prime}(x)=0, \\
\infty, & {\rm otherwise}
\end{array}
\right.\,,\end{equation}

\medskip
{\bf Theorem~B. (\cite{MSS}).} {\it\, Let $n\geqslant 3,$ and let
$\varphi: (0, \infty)\rightarrow [0, \infty)$ be a non-decreasing
Lebesgue measurable function which satisfies Calderon's condition
\begin{equation*}\label{eq7A}\int\limits_{1}^{\infty}\left(\frac{t}{\varphi(t)}\right)^
{\frac{1}{n-2}}dt<\infty\,.
\end{equation*}
Suppose also that there exist constants  $C>0$ and $T>0$ such that
\begin{equation*}\label{eq30}
\varphi(2t)\leqslant C\cdot \varphi(t)\,\,\forall\,\,t\geqslant T\,.
\end{equation*}
Let $Q:{\Bbb B}^n\rightarrow [0, \infty]$ be integrable function in
${\Bbb B}^n.$ Assume that $f$ is a homeomorphism of ${\Bbb B}^n$
onto ${\Bbb B}^n$ such that $f\in W^{1, \varphi}({\Bbb B}^n)$ and,
in addition, $f(0)=0.$
Let, moreover, $K_I(x, f)\leqslant Q(x)$ for a.e. $x\in {\Bbb B}^n$
and, besides that,
\begin{equation*}\label{eq1AE}
\sup\limits_{\varepsilon\in(0,\varepsilon_0)}\frac{1}{\Omega_n\varepsilon^n}\int\limits_{{\Bbb
B}^n\cap B(\zeta,\varepsilon)}Q(x)\,dm(x) < C\qquad \forall\,\,
\zeta \in \partial{\Bbb B}^n
\end{equation*}
holds for some $\varepsilon_0>0,$ where $\Omega_n$ is the volume of
the unit ball in ${\Bbb R}^n.$

Then $f$ has a homeomorphic extension $f:\overline{{\Bbb
B}^n}\rightarrow \overline{{\Bbb B}^n}$ and, in addition,

$$
|f(x_2)-f(x_1)| \leqslant 2\alpha_n\varepsilon_0^{\,-\alpha}\cdot
|x_2-x_1|^{\,\alpha} \qquad \forall\,\, x_1, x_2\in \partial {\Bbb
B}^n: |x_2-x_1|<\delta_0\,,
$$
where $\delta_0:=\min\left\{\frac{1}{2},\varepsilon_0^2\right\}$,
$\omega_{n-1}$ is the area of $n-1$-dimensional sphere
$\mathbb{S}^{n-1}$ and $\alpha:=\left(\frac{\omega_{n-1}\log
2}{\Omega_n(4^n+1)2^{n+1}C}\right)^{1/(n-1)}.$}

\medskip
Let $x_0\in {\Bbb R}^n,$ $0<r_1<r_2<\infty,$
\begin{equation}\label{eq1ED}
S(x_0,r) = \{ x\,\in\,{\Bbb R}^n : |x-x_0|=r\}\,, \quad B(x_0, r)=\{
x\,\in\,{\Bbb R}^n : |x-x_0|<r\}\end{equation}
and
\begin{equation}\label{eq1**}
A=A(x_0, r_1,r_2)=\left\{ x\,\in\,{\Bbb R}^n:
r_1<|x-x_0|<r_2\right\}\,.\end{equation}

\medskip
Recall that a mapping $f:D\rightarrow {\Bbb R}^n$ is called {\it
discrete} if the pre-image $\{f^{-1}\left(y\right)\}$ of each point
$y\,\in\,{\Bbb R}^n$ consists of isolated points, and {\it is open}
if the image of any open set $U\subset D$ is an open set in ${\Bbb
R}^n.$ Later, in the extended space $\overline{{{\Bbb R}}^n}={{\Bbb
R}}^n\cup\{\infty\}$ we use the {\it spherical (chordal) metric}
$h(x,y)=|\pi(x)-\pi(y)|,$ where $\pi$ is a stereographic projection
$\overline{{{\Bbb R}}^n}$ onto the sphere
$S^n(\frac{1}{2}e_{n+1},\frac{1}{2})$ in ${{\Bbb R}}^{n+1},$ namely,
\begin{equation*}\label{eq3C}
h(x,\infty)=\frac{1}{\sqrt{1+{|x|}^2}}\,,\quad
h(x,y)=\frac{|x-y|}{\sqrt{1+{|x|}^2} \sqrt{1+{|y|}^2}}\,, \quad x\ne
\infty\ne y
\end{equation*}
(see \cite[Definition~12.1]{Va}). Further, the closure
$\overline{A}$ and the boundary $\partial A$ of the set $A\subset
\overline{{\Bbb R}^n}$ we understand relative to the chordal metric
$h$ in $\overline{{\Bbb R}^n}.$ We also set
\begin{equation}\label{eq5}
h(E):=\sup\limits_{x,y\in E}h(x, y)\,,\qquad h(A,
B)=\inf\limits_{x\in A, y\in B} h(x, y)\,.
\end{equation}
A Borel function $\rho:{\Bbb R}^n\,\rightarrow [0,\infty] $ is
called {\it admissible} for the family $\Gamma$ of paths $\gamma$ in
${\Bbb R}^n,$ if the relation
\begin{equation*}\label{eq1.4}
\int\limits_{\gamma}\rho (x)\, |dx|\geqslant 1
\end{equation*}
holds for all (locally rectifiable) paths $ \gamma \in \Gamma.$ In
this case, we write: $\rho \in {\rm adm} \,\Gamma .$ Let $p\geqslant
1,$ then {\it $p$-modulus} of $\Gamma $ is defined by the equality
\begin{equation*}\label{eq1.3gl0}
M_p(\Gamma)=\inf\limits_{\rho \in \,{\rm adm}\,\Gamma}
\int\limits_{{\Bbb R}^n} \rho^p (x)\,dm(x)\,.
\end{equation*}

\medskip
We set $M(\Gamma):=M_n(\Gamma).$ Given sets $E,$
$F\subset\overline{{\Bbb R}^n}$ and a domain $D\subset {\Bbb R}^n$
we denote by $\Gamma(E,F,D)$ a family of all paths
$\gamma:[a,b]\rightarrow \overline{{\Bbb R}^n}$ such that
$\gamma(a)\in E,\gamma(b)\in\,F $ and $\gamma(t)\in D$ for $t \in
(a, b).$ A domain $G$ in ${\Bbb R}^n$ is called a {\it quasiextremal
distance domain} (short. $QED$-{\it domain}), if there is a constant
$A_0\geqslant 1,$ such that the inequality
\begin{equation}\label{eq4***}
M(\Gamma(E, F, {\Bbb R}^n))\leqslant A_0\cdot M(\Gamma(E, F, G))
\end{equation}
holds for any continua $E, F\subset G.$

Let $S_i=S(x_0, r_i),$ $i=1,2,$ where spheres $S(x_0, r_i)$ centered
at $x_0$ of the radius $r_i$ are defined in~(\ref{eq1ED}). Let
$Q:{\Bbb R}^n\rightarrow {\Bbb R}^n$ be a Lebesgue measurable
function satisfying the condition $Q(x)\equiv 0$ for $x\in{\Bbb
R}^n\setminus D.$ A mapping $f:D\rightarrow \overline{{\Bbb R}^n}$
is called a {\it ring $Q$-mapping at the point $x_0\in
\overline{D}\setminus \{\infty\}$},\index{ring $Q$-mapping} if the
condition
\begin{equation} \label{eq2*!A}
M(f(\Gamma(S_1, S_2, D)))\leqslant \int\limits_{A\cap D} Q(x)\cdot
\eta^n (|x-x_0|)\, dm(x)
\end{equation}
holds for all $0<r_1<r_2<d_0:=\sup\limits_{x\in D}|x-x_0|$ and all
Lebesgue measurable functions $\eta:(r_1, r_2)\rightarrow [0,
\infty]$ such that
\begin{equation}\label{eq8BC}
\int\limits_{r_1}^{r_2}\eta(r)\,dr\geqslant 1\,.
\end{equation}
A mapping $f$ is called a {\it ring $Q$-mapping in $D,$} if
condition~(\ref{eq2*!A}) is satisfied at every point $x_0\in D,$ and
a {\it ring $Q$-mapping in $\overline{D},$} if the
condition~(\ref{eq2*!A}) holds at every point $x_0\in\overline{D}.$
A mapping $f:D\rightarrow\overline{{\Bbb R}^n}$ between domains
$D\subset{\Bbb R}^n$ and $D^{\,\prime}\subset\overline{{\Bbb R}^n}$
is called {\it closed} if $C(f, \partial D)\subset \partial
D^{\,\prime},$ where, as usual, $C(f, \partial D)$ is the cluster
set of the mapping $f$ on $\partial D.$

\medskip
Given numbers $A_0>0,$ $R_0>0$ and $\delta>0,$ a domain $D\subset
{\Bbb R}^n,$ $n\geqslant 2,$ a point $x_0\in \partial D,$
$x_0\ne\infty,$ and a given function $Q:D\rightarrow[0, \infty],$
$Q(x)\equiv 0$ for $x\in {\Bbb R}^n\setminus D,$ we denote
$\frak{R}^{A_0, R_0}_{Q, \delta}(D, x_0)$ the family of all open,
discrete and closed ring $Q$-mappings $f:D\rightarrow B(0, R_0)$ at
$x_0$ such that the domain $D_f ^{\,\prime}=f(D)$ satisfies the
condition~(\ref{eq4***}) with $G=D_f^{\,\prime},$ and, in addition,
there exists a path connected continuum $K_f\subset D^{\,\prime}_f$
such that ${\rm diam\,}(K_f)\geqslant \delta$ and $h(f^{\,-1}(K_f),
\partial D )\geqslant \delta>0.$ The following theorem holds.

\medskip
{\bf Theorem~C. (\cite{Sev$_4$}).} {\it\, Assume that, the following
conditions hold: 1) there is $r_0=r_0(x_0)>0$ such that, the set
$B(x_0, r)\cap D$ is connected for any $0<r<r_0;$ 2) there is
$0<C=C(x_0)<\infty,$ such that
\begin{equation*}\label{eq1E}
\limsup\limits_{\varepsilon\rightarrow
0}\frac{1}{\Omega_n\cdot\varepsilon^n}\int\limits_{B(x_0,
\varepsilon)\cap D}Q(x)\,dm(x)\leqslant C\,.
\end{equation*}
Then there is
$\widetilde{\varepsilon_0}=\widetilde{\varepsilon}_0(x_0)>0$ and a
number $\alpha=\alpha(n, \delta, C, R_0, x_0)>0$ such that the
relation
\begin{equation*}\label{eq2}
|f(x)-f(y)|\leqslant \alpha\cdot\max\{|x-x_0|^{\beta},
|y-x_0|^{\beta}\}\end{equation*}
holds for any $x, y\in B(x_0, \widetilde{\varepsilon}(x_0))\cap D$
and all $f\in\frak{R}^{A_0, R_0}_{Q, \delta}(D, x_0),$ where
$\beta=\left(\frac{n\log 2}{A_0C2^{n+1}}\right)^{\frac{1}{n-1}}.$}

\medskip
Note that in Theorem~\textbf{A}, the mapping~$f$ is quasiregular, in
Theorem~\textbf{B} it is a homeomorphism, and in Theorem~\textbf{C}
it is closed (it may have branch points, but preserves the boundary
of the domain). {\bf Our further goal is to abandon the
quasiregularity condition, the homeomorphism condition, and also the
preservation of the boundary of a domain.} That is, we should
consider a case that is more general in comparison with all those
described ``in a certain sense''; this is what the entire further
presentation is devoted to. The key point is to avoid the boundary
preservation under mapping. We emphasize that it is impossible to
completely avoid the preservation of the boundary, since even for
isolated boundary points, the mapping does not always have a
continuous extension. Nevertheless, it is possible to identify broad
classes of mappings and domains for which H\"{o}lder continuity or
more general distortion estimates remain valid even if the boundary
of a domain is not preserved. In the statement below we will
consider one of such classes.

\medskip
Let
\begin{gather}\nonumber
N(y, f, A)\,=\,{\rm card}\,\left\{x\in A: f(x)=y\right\}\,,\\
\label{eq1_A_3} N(f, A)\,=\,\sup\limits_{y\in{\Bbb R}^n}\,N(y, f,
A)\,.
\end{gather}
Given $N\in {\Bbb N},$ $\delta>0,$ $A_0>0,$ sets $E, E_*$ in
$\overline{{\Bbb R}^n},$ $n\geqslant 3,$ while $E_*$ is closed in
$\overline{{\Bbb R}^n},$ a domain $D\subset {\Bbb R}^n,$ an
increasing function $\varphi:[0,\infty)\rightarrow[0,\infty)$ and a
Lebesgue measurable function $Q:D\rightarrow [0, \infty]$ we denote
by $\frak{R}^{E_*, E, N}_{\varphi, Q, \delta, A_0}(D)$ the family of
all bounded open discrete mappings $f:D\rightarrow {\Bbb R}^n$ in
the class $W_{\rm loc}^{1, \varphi}(D)$ such that $K_I(x,
f)\leqslant Q(x)$ a.e., $N(f, D)\leqslant N$ and, in addition,

\medskip
1) $C(f, \partial D)\subset E_*,$

\medskip
2) for each component $K$ of $D^{\,\prime}_f\setminus  E_*,$
$D^{\,\prime}_f:=f(D),$ there is a continuum $K_f\subset K$ such
that $h(K_f)\geqslant \delta$ and $h(f^{\,-1}(K_f), \partial
D)\geqslant \delta>0,$

\medskip
3) $f^{\,-1}(E_*)=E,$

\medskip
4) any component $G$ of $f(D)\setminus E_*$ satisfies the
relation~(\ref{eq4***}).

\medskip
Given $N\in {\Bbb N},$ $\delta>0,$ $A_0>0,$ sets $E, E_*$ in
$\overline{{\Bbb R}^2},$ while $E_*$ is closed in $\overline{{\Bbb
R}^2},$ a domain $D\subset {\Bbb R}^2$ and a Lebesgue measurable
function $Q:D\rightarrow [0, \infty]$ we denote by $\frak{R}^{E_*,
E, N}_{Q, \delta, A_0}(D)$ the family of all bounded open discrete
mappings $f:D\rightarrow {\Bbb R}^2$ in the class $W_{\rm loc}^{1,
1}(D)$ such that $K_I(x, f)\leqslant Q(x)$ a.e., $N(f, D)\leqslant
N$ and, in addition, the relations 1)--4) mentioned above hold. The
following statement is true.

\medskip
\begin{theorem}\label{th4} {\it\, Let
$D$ be a bounded domain in ${\Bbb R}^n,$ $n\geqslant 2,$ $x_0\in
\partial D.$ Assume that:

\medskip
1) the set $E$ is nowhere dense in $D,$ and $D$ is finitely
connected on $E\cup \partial D,$ i.e., for any $z_0\in E\cup
\partial D$ and any neighborhood $\widetilde{U}$ of $z_0$ there is a
neighborhood $\widetilde{V}\subset \widetilde{U}$ of $z_0$ such that
$(D\cap \widetilde{V})\setminus E$ consists of finite number of
components;

\medskip
2) there are $m=m(x_0)\in {\Bbb N},$ $1\leqslant m<\infty,$ and
$r_0=r_0(x_0)>0$ such that the following is true:

\medskip
2a) $B(x_0, r)\cap D$ is connected for any $0<r<r_0;$

\medskip
2b) $(B(x_0, r)\cap D)\setminus E$ consists at most of $m$
components for any $0<r<r_0;$ besides that, for any component $K_1$
in $(B(x_0, r)\cap D)\setminus E$ and any $x, y\in K_1$ there is a
path $\gamma:[a, b]\rightarrow {\Bbb R}^n$ such that $|\gamma|\in
K_1\cap \overline{B(x_0 , \max\{|x-x_0|, |y-x_0|\})};$

\medskip
3) there is $0<C=C(x_0)<\infty$ such that
\begin{equation}\label{eq1E_1}
\limsup\limits_{\varepsilon\rightarrow
0}\frac{1}{\Omega_n\cdot\varepsilon^n}\int\limits_{B(x_0,
\varepsilon)\cap D}Q(x)\,dm(x)\leqslant C\,;
\end{equation}
4) if $n\geqslant 3,$ assume that, the function $\varphi$ satisfies
Calderon condition
\begin{equation}\label{eq1_A_10}
\int\limits_{1}^{\infty}\left(\frac{t}{\varphi(t)}\right)^
{\frac{1}{n-2}}\,dt<\infty\,.
\end{equation}
Then there are $\widetilde{\varepsilon}_0>0$ and $\alpha>0$ such
that the relation
\begin{equation*}\label{eq2_3}
h(f(x), f(y))\leqslant \alpha\cdot\max\{|x-x_0|^{\beta},
|y-x_0|^{\beta}\}\end{equation*}
holds for any $x, y\in B(x_0, \widetilde{\varepsilon}(x_0))\cap D$
and all $f\in\frak{R}^{E_*, E, N}_{\varphi, Q, \delta, A_0}(D)$
whenever $n\geqslant 3$ (all $f\in\frak{R}^{E_*, E, N}_{Q, \delta,
A_0}(D)$ whenever $n=2$), where $0<\beta\leqslant 1$ is some
constant depending only on $n,$ $N,$ $C$ and $A_0.$ }
\end{theorem}

\medskip
The statement similar to Theorem~\ref{th3} holds for ``domains with
prime ends'', see below. The necessary definitions related to ends,
prime ends, chains of cuts, etc. can be found, for example, in
Chapter~5 of monograph~\cite{Sev$_3$}. Let $K$ be the end of $D$ in
${\Bbb R}^n$, then the set
\begin{equation}\label{eq1A_7}
I(K)=\bigcap\limits_{m=1}\limits^{\infty}\overline{d_m}
\end{equation}
is called {\it the impression of the end} $K$.
Following~\cite{Na$_2$}, we say that the end $K$ is {\it a prime
end}, if $K$ contains a chain of cuts $\{\sigma_m\}$ such that
$$\lim\limits_{m\rightarrow\infty}M(\Gamma(C, \sigma_m, D))=0$$
for some continuum $C$ in $D.$ In the following, the following
notation is used: the set of prime ends corresponding to the domain
$D,$ is denoted by $E_D,$ and the completion of the domain $D$ by
its prime ends is denoted $\overline{D}_P.$

\medskip
Consider the following definition, which goes back to
N\"akki~\cite{Na$_2$}, cf.~\cite{KR}. The boundary of a domain $D$
in ${\Bbb R}^n$ is said to be {\it locally quasiconformal}if every
$x_0\in\partial D$ has a neighborhood $U$ that admits a
quasiconformal mapping $\kappa$ onto the unit ball ${\Bbb
B}^n\subset{\Bbb R}^n$ such that $\kappa(\partial D\cap U)$ is the
intersection of ${\Bbb B}^n$ and a coordinate hyperplane. The
sequence of cuts $\sigma_m,$ $m=1,2,\ldots ,$ is called {\it
regular,} if
$\overline{\sigma_m}\cap\overline{\sigma_{m+1}}=\varnothing$ for
$m\in {\Bbb N}$ and, in addition, $d(\sigma_{m})\rightarrow 0$ as
$m\rightarrow\infty.$ If the end $K$ contains at least one regular
chain, then $K$ will be called {\it regular}. We say that a bounded
domain $D$ in ${\Bbb R}^n$ is {\it regular}, if $D$ can be
quasiconformally mapped to a domain with a locally quasiconformal
boundary whose closure is a compact in ${\Bbb R}^n,$ and, besides
that, every prime end in $D$ is regular. Note that space
$\overline{D}_P=D\cup E_D$ is metric, which can be demonstrated as
follows. If $g:D_0\rightarrow D$ is a quasiconformal mapping of a
domain $D_0$ with a locally quasiconformal boundary onto some domain
$D,$ then for $x, y\in \overline{D}_P$ we put:
\begin{equation}\label{eq5M}
\rho(x, y):=|g^{\,-1}(x)-g^{\,-1}(y)|\,,
\end{equation}
where the element $g^{\,-1}(x),$ $x\in E_D,$ is to be understood as
some (single) boundary point of the domain $D_0.$ The specified
boundary point is unique and well-defined by~\cite[Theorem~2.1,
Remark~2.1]{IS}, cf.~\cite[Theorem~4.1]{Na$_2$}. It is easy to
verify that~$\rho$ in~(\ref{eq5M}) is a metric on $\overline{D}_P.$
If $g_*$ is another quasiconformal mapping of a domain $D_*$ with
locally quasiconformal boundary onto $D$, then the corresponding
metric
$\rho_*(p_1,p_2)=|{\widetilde{g_*}}^{-1}(p_1)-{\widetilde{g_*}}^{-1}(p_2)|$
generates the same convergence and, consequently, the same topology
in $\overline {D}_P$ as $\rho_0$ because $g_0\circ g_*^{-1}$ is a
quasiconformal mapping of $D_*$ and $D_0$, which extends, by Theorem
4.1 in~\cite{Na$_2$}, to a homeomorphism between $\overline {D_*}$
and $\overline {D_0}$. In the sequel, this topology in $\overline
{D}_P$ will be called the {\it topology of prime ends}; the
continuity of mappings $F\colon \overline
{D}_P\rightarrow\overline{D^{\,\prime}}_P$ will be understood
relative to this topology. The following statement holds.

\medskip
\begin{theorem}\label{th5} {\it\, Let
$D$ be a regular domain in ${\Bbb R}^n,$ $n\geqslant 2,$ $P_0\in
E_D.$ Assume that:

\medskip
1) the set $E$ is nowhere dense in $D,$ and $D$ is finitely
connected on $E,$ i.e., for any $z_0\in E$ and any neighborhood
$\widetilde{U}$ of $z_0$ there is a neighborhood
$\widetilde{V}\subset \widetilde{U}$ of $z_0$ such that $(D\cap
\widetilde{V})\setminus E$ consists of finite number of components;

2) for each $y_0\in\partial D$ there exists $r_0=r_0(y_0)
>0$ such that:

2a) there is $p=p(y_0)\in {\Bbb N}$ such that the set $B(y_0, r)\cap
D$ consists at most of $p$ components for all $0<r<r_0,$ and, for
each component $K$ of the set $B(y_0, r)\cap D$ the following
condition is fulfilled: any $x, y\in K$ may be joined by a path
$\gamma:[a, b]\rightarrow {\Bbb R}^n$ such that $|\gamma|\in K\cap
\overline{B(y_0 , \max\{|x-y_0|, |y-y_0|\})},$
$$|\gamma|=\{x\in {\Bbb R}^n: \exists\,t\in[a, b]: \gamma(t)=x\};$$

2b) there is $m=m(y_0)\in {\Bbb N}$ such that, for $K$ mentioned
above, $K\setminus E$ consists at most of $m$ components for any
$0<r<r_0$ and, for each component $K_1$ of the set $K\setminus E$
the following condition is fulfilled: any $x, y\in K_1$ may be
joined by a path $\gamma:[a, b]\rightarrow {\Bbb R}^n$ such that
$|\gamma|\in K_1\cap \overline{B(y_0 , \max\{|x-y_0|, |y-y_0|\})};$

\medskip
3) for any $x_0\in I(P_0)$ there is $0<C=C(x_0)<\infty$ such that
\begin{equation}\label{eq1EA}
\limsup\limits_{\varepsilon\rightarrow
0}\frac{1}{\Omega_n\cdot\varepsilon^n}\int\limits_{B(x_0,
\varepsilon)\cap D}Q(x)\,dm(x)\leqslant C\,,
\end{equation}
where $I(P_0)$ denotes the impression of $P_0,$ see~(\ref{eq1A_7}).

4) If $n\geqslant 3,$ assume that, the function $\varphi$ satisfies
Calderon condition~(\ref{eq1_A_10}).

Then for each $P\in E_D$ there exists $y_0\in \partial D$ such that
$I(P)=y_0,$ where $I(P)$ denotes the impression of $P$ defined
in~(\ref{eq1A_7}). In addition, there exists a neighborhood $U$ of
$P_0$ in the metric space $(\overline{D}_P, \rho)$ such that the
inequality
\begin{equation*}\label{eq2_1}
h(f(x), f(y))\leqslant \alpha\cdot\max\{|x-x_0|^{\beta},
(|y-x_0|)^{\beta}\}\end{equation*}
holds for any $x, y\in U\cap D$ and all $f\in\frak{R}^{E_*, E,
N}_{\varphi, Q, \delta, A_0}(D)$ whenever $n\geqslant 3$ (all
$f\in\frak{R}^{E_*, E, N}_{Q, \delta, A_0}(D)$ whenever $n=2$),
where $x_0:=I(P_0)$ and $0<\beta\leqslant 1$ is some constant
depending only on $n,$ $N,$ $C$ and $A_0.$}
\end{theorem}

\section{Preliminaries}

Below, we assume that the definitions of the moduli of families of
paths and surfaces used below are known. All these definitions, in
fairly detailed form, can be found in~\cite{Sev$_3$}. The next class
of mappings is a generalization of quasiconformal mappings in the
sense of Gehring's ring definition (see \cite{Ge}; it is the subject
of a separate study, see, e.g., \cite[Chapter~9]{MRSY}). Let $D$ and
$D^{\,\prime}$ be domains in ${\Bbb R}^n$ with $n\geqslant 2$.
Suppose that $x_0\in\overline {D}\setminus\{\infty\}$ and $Q\colon
D\rightarrow(0,\infty)$ is a Lebesgue measurable function. A mapping
$f:D\rightarrow D^{\,\prime}$ is called a {\it lower $Q$-mapping at
a point $x_0$ relative to the $p$-modulus} if
\begin{equation}\label{eq1A}
M_p(f(\Sigma_{\varepsilon}))\geqslant \inf\limits_{\rho\in{\rm
ext}_p\,{\rm adm}\Sigma_{\varepsilon}}\int\limits_{D\cap A(x_0,
\varepsilon, r_0)}\frac{\rho^p(x)}{Q(x)}\,dm(x)
\end{equation}
for every spherical ring $A(x_0, \varepsilon, r_0)=\{x\in {\Bbb
R}^n\,:\, \varepsilon<|x-x_0|<r_0\}$, $r_0\in(0, d_0)$,
$d_0=\sup_{x\in D}|x-x_0|$, where $\Sigma_{\varepsilon}$ is the
family of all intersections of the spheres $S(x_0, r)$ with the
domain $D$, $r\in (\varepsilon, r_0)$. If $p=n$, we say that $f$ is
a lower $Q$-mapping at $x_0$. We say that $f$ is a lower $Q$-mapping
relative to the $p$-modulus in $A\subset \overline {D}$ if
(\ref{eq1A}) is true for all $x_0\in A$. The following lemma holds,
see e.g. \cite[Lemma~5.1]{Sev$_3$}.
\begin{lemma}\label{thOS4.1} {\it\,Let $D$ be a domain in ${\Bbb R}^n$,
$n\geqslant  2$, and let $\varphi\colon (0,\infty)\rightarrow
(0,\infty)$ be a monotone nondecreasing function satisfying
(\ref{eq1_A_10}). If $p>n-1$, then every open discrete mapping
$f\colon D\rightarrow {\Bbb R}^n$ of class $W^{1,\varphi}_{{\rm
loc}}$ with finite distortion and such that $N(f, D)<\infty$ is a
lower $Q$-mapping relative to the $p$-modulus at every point
$x_0\in\overline {D}$ for
$$
Q(x)=N(f, D)\cdot K^{\frac{p-n+1}{n-1}}_{I, \alpha}(x, f),
$$
$\alpha:=\frac{p}{p-n+1}$, where the inner dilation $K_{I,\alpha}(x,
f)$ for~$f$ at~$x$ of order~$\alpha$ is defined by \eqref{eq0.1.1A},
and the multiplicity $N(f, D)$ is defined by the
relation~(\ref{eq1_A_3}).}
\end{lemma}

\medskip
The statement similar to Lemma~\ref{thOS4.1} holds for $n=2,$ but
for Sobolev classes (see e.g. \cite[Theorem~4]{Sev$_2$}).

\medskip
\begin{lemma}\label{lem3A} {\it\,Let $D$ be a domain in ${\Bbb R}^2$,
$n\geqslant  2$, and let $p>1.$ Then every open discrete mapping
$f\colon D\rightarrow {\Bbb R}^2$ of class $W^{1,1}_{{\rm loc}}$
with finite distortion and such that $N(f, D)<\infty$ is a lower
$Q$-mapping relative to the $p$-modulus at every point
$x_0\in\overline {D}$ for
$$
Q(x)=N(f, D)\cdot K^{p-1}_{I, \alpha}(x, f),
$$
$\alpha:=\frac{p}{p-1}$, where the inner dilation $K_{I,\alpha}(x,
f)$ for~$f$ at~$x$ of order~$\alpha$ is defined by \eqref{eq0.1.1A},
and the multiplicity $N(f, D)$ is defined by the
relation~(\ref{eq1_A_3}).}
\end{lemma}

\medskip
Let $S_i=S(x_0, r_i),$ $i=1,2,$ where spheres $S(x_0, r_i)$ centered
at $x_0$ of the radius $r_i$ are defined in~(\ref{eq1ED}). The
following statement holds, see~\cite{Sev$_1$}.

\medskip
\begin{lemma}\label{lem4}
{\,\it Let $x_0\in \partial D,$ let $f:D\rightarrow {\Bbb R}^n$ be a
bounded, open, discrete, and closed lower $Q$-mapping with respect
to $p$-modulus in domain $D\subset{\Bbb R}^n,$ $Q\in
L_{loc}^{\frac{n-1}{p-n+1}}({\Bbb R}^n),$ $n-1<p,$ and
$\alpha:=\frac{p}{p-n+1}.$ Then, for every
$\varepsilon_0<d_0:=\sup\limits_{x\in D}|x-x_0|$ and every compact
set $C_2\subset D\setminus B(x_0, \varepsilon_0)$ there exists
$\varepsilon_1,$ $0<\varepsilon_1<\varepsilon_0,$ such that, for
each $\varepsilon\in (0, \varepsilon_1)$ and each compactum
$C_1\subset \overline{B(x_0, \varepsilon)}\cap D$ the inequality
\begin{equation*}\label{eq3A_2}
M_{\alpha}(f(\Gamma(C_1, C_2, D)))\leqslant \int\limits_{A(x_0,
\varepsilon, \varepsilon_1)}Q^{\frac{n-1}{p-n+1}}(x)
\eta^{\alpha}(|x-x_0|)\,dm(x)
\end{equation*}
holds, where $A(x_0, \varepsilon, \varepsilon_1)=\{x\in {\Bbb R}^n:
\varepsilon<|x-x_0|<\varepsilon_1\}$ and $\eta: (\varepsilon,
\varepsilon_1)\rightarrow [0,\infty]$ is an arbitrary Lebesgue
measurable function such that
\begin{equation}\label{eq6B}
\int\limits_{\varepsilon}^{\varepsilon_1}\eta(r)\,dr=1\,.
\end{equation}
}
\end{lemma}
Let $D\subset {\Bbb R}^n,$ $f:D\rightarrow {\Bbb R}^n$ be a discrete
open mapping, $\beta: [a,\,b)\rightarrow {\Bbb R}^n$ be a path, and
$x\in\,f^{\,-1}(\beta(a)).$ A path $\alpha: [a,\,c)\rightarrow D$ is
called a {\it maximal $f$-lifting} of $\beta$ starting at $x,$ if
$(1)\quad \alpha(a)=x\,;$ $(2)\quad f\circ\alpha=\beta|_{[a,\,c)};$
$(3)$\quad for $c<c^{\prime}\leqslant b,$ there is no a path
$\alpha^{\prime}: [a,\,c^{\prime})\rightarrow D$ such that
$\alpha=\alpha^{\prime}|_{[a,\,c)}$ and $f\circ
\alpha^{\,\prime}=\beta|_{[a,\,c^{\prime})}.$ Here and in the
following we say that a path $\beta:[a, b)\rightarrow
\overline{{\Bbb R}^n}$ converges to the set $C\subset
\overline{{\Bbb R}^n}$ as $t\rightarrow b,$ if $h(\beta(t),
C)=\sup\limits_{x\in C}h(\beta(t), C)\rightarrow 0$ at $t\rightarrow
b.$ The following is true (see~\cite[Lemma~3.12]{MRV$_2$}).

\medskip
\begin{proposition}\label{pr3}
{\it\, Let $f:D\rightarrow {\Bbb R}^n,$ $n\geqslant 2,$ be an open
discrete mapping, let $x_0\in D,$ and let $\beta: [a,\,b)\rightarrow
{\Bbb R}^n$ be a path such that $\beta(a)=f(x_0)$ and such that
either $\lim\limits_{t\rightarrow b}\beta(t)$ exists, or
$\beta(t)\rightarrow \partial f(D)$ as $t\rightarrow b.$ Then
$\beta$ has a maximal $f$-lifting $\alpha: [a,\,c)\rightarrow D$
starting at $x_0.$ If $\alpha(t)\rightarrow x_1\in D$ as
$t\rightarrow c,$ then $c=b$ and $f(x_1)=\lim\limits_{t\rightarrow
b}\beta(t).$ Otherwise $\alpha(t)\rightarrow \partial D$ as
$t\rightarrow c.$}
\end{proposition}

\medskip
The following statement holds, see, e.g., \cite[Theorem~1]{D},
cf.~\cite[Lemma~2.1]{DS}.

\medskip
\begin{lemma}\label{lem2}
{\it\, Let $D$ be a domain in ${\Bbb R}^n,$ $n\geqslant 2,$ and let
$x_0\in \partial D.$ Assume that $E$ is closed and nowhere dense in
$D,$ and $D$ is finitely connected on $E,$ i.e., for any $z_0\in E$
and any neighborhood $\widetilde{U}$ of $z_0$ there is a
neighborhood $\widetilde{V}\subset \widetilde{U}$ of $z_0$ such that
$(D\cap \widetilde{V})\setminus E$ consists of finite number of
components.

In addition, assume that the following condition holds: for any
$x_0\in\partial D$ there are $m=m(x_0)\in {\Bbb N},$ $1\leqslant
m<\infty,$ $r_0=r_0(x_0)>0$ and a neighborhood $U$ of $x_0$ such
that the following is true: $B(x_0, r_0)\subset U$ and

\medskip
2a) $B(x_0, r)\cap D$ is connected for any $0<r<r_0;$

\medskip
2b) $(B(x_0, r)\cap D)\setminus E$ consists at most of $m$
components;

Let $x_k, y_k\in D\setminus E,$ $k=1,2,\ldots ,$ be a sequences
converging to $x_0$ as $k\rightarrow\infty.$ Then there are paths
$\gamma_k,$ $k=1,2,\ldots ,$ such that $\gamma_k$ lies in
$\overline{B(x_0, r_k)}\cap D,$ $r_k=\max\{|x_k-x_0|, |y_k-x_0|\},$
$\gamma_k:[0, 1]\rightarrow \overline{B(x_0, r_k)}\cap D,$
$\gamma_k(0)=x_k,$ $\gamma_k(1)=y_k,$ $\gamma_k(t)\in B(x_0,
r_k)\cap D$ for $0<t<1,$ whenever $|\gamma_k|:=\{x\in {\Bbb R}^n:
\exists\,t\in [0, 1]: \gamma_k(t)=x\}$ contains at most $m-1$ points
in~$E.$ }
\end{lemma}

\medskip
A domain $R$ in $\overline{{\Bbb R}^n},$ $n\geqslant 2,$ is called
{\it a ring,} if $\overline{{\Bbb R}^n}\setminus R$ consists of
exactly two components $E$ and $F.$ In this case, we write: $R=R(E,
F).$ The following statement is true, see \cite[ratio~(7.29)]{MRSY}.

\medskip
\begin{proposition}\label{pr3_1}{\it\, If $R=R(E, F)$ is a ring, then
$$M(\Gamma(E, F, \overline{{\Bbb R}^n}))\geqslant\frac{\omega_{n-1}}{\left(
\log\frac{2\lambda^2_n}{h(E)h(F)}\right)^{n-1}}\,,$$
where $\lambda_n \in[4,2e^{n-1}),$ $\lambda_2=4$ and
$\lambda_n^{1/n} \rightarrow e$ as $n\rightarrow \infty,$ and $h(E)$
denotes the chordal diameter of the set $E$ defined in~(\ref{eq5}).}
\end{proposition}

\medskip
Let $a>0$ and let $\varphi\colon [a,\infty)\rightarrow[0,\infty)$ be
a nondecreasing function such that, for some constants $\gamma>0,$
$T>0 $ and all $t\geqslant T$, the inequality
\begin{equation}\label{eq1B}
\varphi(2t)\leqslant \gamma\cdot\varphi(t)
\end{equation}
is fulfilled. We will call such functions {\it functions that
satisfy the doubling condition}.

Let $\varphi\colon [a,\infty)\rightarrow[0,\infty)$ be a function
with the doubling condition, then the function
$\widetilde{\varphi}(t):=\varphi(1/t)$ does not increase and is
defined on a half-interval $(0, 1/a].$ The following statement is
proved in~\cite[Lemma~3.1]{RSS}.

\medskip
\begin{proposition}\label{pr1}{\it\,
Let $a>0,$ let $\widetilde{\varphi}\colon
[a,\infty)\rightarrow[0,\infty)$ be a nondecreasing function with a
doubling condition~(\ref{eq1B}), let $x_0\in {\Bbb R}^n,$
$n\geqslant 2,$ and let $Q:{\Bbb R}^n\rightarrow [0, \infty]$ be a
Lebesgue measurable function for which there exists $0<C<\infty$
such that
\begin{equation*}\label{eq1AA}
\limsup\limits_{\varepsilon\rightarrow
0}\frac{\widetilde{\varphi}(1/\varepsilon)}{\Omega_n\cdot\varepsilon^n
}\int\limits_{B(x_0, \varepsilon)}Q(x)\,dm(x)\leqslant C\,.
\end{equation*}
Then there exists $\varepsilon^{\,\prime}_0>0$ such that
\begin{equation*}\label{eq2B_2}
\int\limits_{\varepsilon<|x-x_0|<\varepsilon^{\,\prime}_0}
\frac{\widetilde{\varphi}(1/|x-x_0|)Q(x)\,dm(x)}{|x-x_0|^n}\leqslant
C_1\cdot\left(\log\frac{1}{\varepsilon}\right)\,,\qquad\varepsilon\rightarrow
0\,,
\end{equation*}
where $C_1:=\frac{\gamma C\Omega_n2^n}{\log 2}.$}
\end{proposition}

\section{Main Lemmas}

The following statement holds.

\medskip
\begin{lemma}\label{lem1}
{\it\, Let $D$ be a domain in ${\Bbb R}^n,$ $n\geqslant 2,$ let $E,
E_*$ be sets in $\overline{{\Bbb R}^n},$ while $E_*$ is closed in
$\overline{{\Bbb R}^n}$ and let $f:D\rightarrow {\Bbb R}^n$ be an
open discrete mapping such that $C(f, \partial D)\subset E_*$ and
$f^{\,-1}(E_*)=E.$ Assume that, for any $x_0\in\partial D$ there are
$m=m(x_0)\in {\Bbb N},$ $1\leqslant m<\infty,$ $r_0=r_0(x_0)>0$ such
that $(B(x_0, r)\cap D)\setminus E$ consists at most of $m$
components for any $0<r<r_0.$

\medskip
Then:

\medskip
1) $D\setminus E$ consists of finite number of components
$D_1,\ldots, D_s,$ $1\leqslant s<\infty,$

\medskip
2) $f$ is closed in each $D_i,$ $i=1,2,\ldots, $ i.e. $f(S)$ is
closed in $f(D_i)$ whenever $S$ is closed in $D_i,$

\medskip
3)  If $f(D_i)\subset K,$ where $K$ is a component of $f(D)\setminus
E_*,$ then $f(D_i)=K.$}
\end{lemma}

\medskip
\begin{proof} a) Observe that, $E$ is closed in $D$ due to the
relation $f^{\,-1}(E_*)=E,$ because $E_*$ is closed in
$\overline{{\Bbb R}^n}$ by the assumption, see e.g.
\cite[Theorem~1.IV.13.1]{Ku$_1$}. Let us to prove that $D\setminus
E$ consists of finite number of components $D_1,\ldots, D_s,$
$1\leqslant s<\infty.$ Indeed, assume the contrary, namely, that
$D\setminus E$ consists of infinite number of components $D_1,
D_2,\ldots .$ Let $z_i\in D_i,$ $i=1,2,\ldots .$ Since
$\overline{{\Bbb R}^n}$ is a compact space, there is a subsequence
$z_{i_k},$ $k=1,2,\ldots ,$ which converge to some point $z_0\in
\overline{D}.$ By the assumption, there is a neighborhood $V=B(z_0,
r),$ $r>0,$ of the point $z_0$ such that $(V\cap D)\setminus E$
consists of no more than $m$ components. Thus, there is at least one
such a component $V_{m_0}$ intersecting infinitely many $D_i,$
contradiction. Thus, $D\setminus E=D_1,\ldots, D_s,$ $1\leqslant
s<\infty,$ as required.

\medskip
b) Now, let us to prove that $f$ is closed in each $D_i,$
$i=1,2,\ldots, $ i.e. $f(S)$ is closed in $f(D_i)$ whenever $S$ is
closed in $D_i.$ Since $f$ is open and discrete, it is sufficient to
prove that $f$ is boundary preserving, that is $C(f,
\partial D_i)\subset \partial f(D_i)$ (see, e.g.,
\cite[Theorem~3.3]{Vu}). Indeed, let $z_k\in D_i,$ $k=1,2,\ldots ,$
let $z_k\rightarrow z_0\in\partial D_i$ and let $f(z_k)\rightarrow
w_0$ as $k\rightarrow\infty.$ We need to prove that $w_0\in
\partial f(D_i).$ Assume the contrary, i.e. $w_0\in f(D_i).$
There are two cases: $z_0\in \partial D$ or $z_0\in D.$ 1) In the
first case, when $z_0\in \partial D,$ we obtain that $w_0\in C(f,
\partial D)\subset E_*.$ But since $w_0\in f(D_i),$ there is
$\zeta_i\in D_i$ such that $f(\zeta_i)=w_0.$ Now, since $w_0\in
E_*,$ we obtain that $\zeta_i\in f^{\,-1}(E_*)=E$ and, consequently,
$\zeta_i\not\in D_i,$ contradiction. 2) Let us consider the second
case, when $z_0\in D.$ Now, by the definition of $D_i,$ we have that
$z_0\in E.$ Therefore, $w_0\in f(E)\subset E_*.$ But since $w_0\in
f(D_i),$ there is $\zeta_i\in D_i$ such that $f(\zeta_i)=w_0.$ Since
$w_0\in E_*,$ we have that $\zeta_i\in f^{\,-1}(E_*)=E$ and,
consequently, $\zeta_i\not\in D_i,$ contradiction.

\medskip
c) Let $D_i$ be such that $f(D_i)\subset K,$ where $K$ is a
component of $f(D)\setminus E_*.$ Observe that $f(D_i)=K.$ Indeed,
$f(D_i)\subset K$ by the definition. Let us to prove that $K\subset
f(D_i).$ Let us prove this inclusion by the contradiction, i.e., let
$a_0\in K\setminus f(D_i).$ Chose $b_0\in f(D_i)$ and join the
points $b_0$ and $a_0$ by a path $\beta:[0, 1]\rightarrow K.$ Let
$\alpha:[0, c)\rightarrow D$ be a maximal $f$-lifting of $\beta$
starting at $c_0:=f^{\,-1}(b_0)\cap D_i$ (this lifting exists by
Proposition~\ref{pr3}). By the same proposition either one of two
cases are possible: $\alpha(t)\rightarrow x_1$ as $t\rightarrow c,$
or $\alpha(t)\rightarrow \partial D_i$ as $t\rightarrow c.$ In the
first case, by Proposition~\ref{pr3} we obtain that $c=1$ and
$f(\beta(1))=f(x_1)=a_0$ which contradict with the choice of $a_0.$
In the second case, when $\alpha(t)\rightarrow \partial D_i$ as
$t\rightarrow c,$ we obtain that $f(\beta(c))\in C(f, \partial
D_i)\subset E_*.$ The latter contradicts the definition of $\beta$
because $\beta$ does not contain itself points in $E_*.$
\end{proof}

\medskip
The following lemma contains the statement about H\"{o}lder
continuity of mappings for ``good'' boundaries. Similar statements
were proved for homeomorphisms and open discrete closed mappings,
see e.g. \cite{MSS}, \cite{RSS} and \cite{Sev$_4$}.

\medskip
\begin{lemma}\label{lem3} {\it\, Assume
that, under conditions of Theorem~\ref{th4}, instead of the
relation~(\ref{eq1E_1}) the following is true: there are
$\varepsilon_0>0,$ $1\leqslant p<n,$ $K_0>0,$ and a Lebesgue
measurable function $\psi:(0, \varepsilon_0)\rightarrow [0, \infty]$
such that
\begin{equation}\label{eq3AB} 0<I(\varepsilon, \varepsilon_0):=
\int\limits_{\varepsilon}^{\varepsilon_0}\psi(t)\,dt < \infty\,,
\qquad I(\varepsilon, \varepsilon_0)\rightarrow \infty\,,\quad
\varepsilon\rightarrow 0\,,\end{equation}
while
\begin{equation} \label{eq3.7A}
\int\limits_{\varepsilon<|x-x_0|<\varepsilon_0}
Q(x)\cdot\psi^n(|x-x_0|)\,dm(x) \leqslant K_0\cdot I^p(\varepsilon,
\varepsilon_0)
\end{equation}
as $\varepsilon\rightarrow 0.$
Then there are $\widetilde{\varepsilon}_0>0$ and $\alpha>0$ such
that the relation
\begin{equation}\label{eq2_2}
h(f(x), f(y))\leqslant \alpha\cdot\max\{\theta(|x-x_0|),
\theta(|y-x_0|)\}\,,\end{equation}
$$\theta(\varepsilon):=\exp\{-(A_0K_0N\omega^{\,-1}_{n-1})^{-\frac{1}{n-1}}\cdot
I^{\frac{n-p}{n-1}}(\varepsilon, \varepsilon_0)\}\,,$$
holds for any $x, y\in B(x_0, \widetilde{\varepsilon}(x_0)\cap D$
and all $f\in\frak{R}^{E_*, E, N}_{\varphi, Q, \delta, A_0}(D)$
whenever $n\geqslant 3$ (all $f\in\frak{R}^{E_*, E, N}_{Q, \delta,
A_0}(D)$ whenever $n=2$).}
\end{lemma}

\medskip
\begin{proof}
{\bf I.} Observe that, $E$ is closed in $D$ due to the relation
$f^{\,-1}(E_*)=E,$ because $E_*$ is closed in $\overline{{\Bbb
R}^n}$ by the assumption, see e.g. \cite[Theorem~1.IV.13.1]{Ku$_1$}.
It is sufficiently to prove Lemma~\ref{lem3} for $x, y\in D\setminus
E,$ because the set $E$ is nowhere dense in $D,$ so that we may
obtain~(\ref{eq2_2}) by means of limit transition as
$\widetilde{x}\rightarrow x$ and $\widetilde{y}\rightarrow y,$ where
$\widetilde{x}, \widetilde{y}\in D\setminus E.$ Let us prove by
contradiction, i.e., assume that the statement of Theorem~\ref{th4}
is not true for $x, y\in D\setminus E.$ Now, there are sequences
$x_k, y_k\rightarrow x_0$ as $k\rightarrow\infty,$ $x_k, y_k\in
D\setminus E,$ and $f_k\in \frak{R}^{E_*, E, N}_{\varphi, Q, \delta,
A_0}(D)$ whenever $n\geqslant 3$ ($f_k\in\frak{R}^{E_*, E, N}_{Q,
\delta, A_0}(D)$ whenever $n=2$), such that
\begin{equation}\label{eq2A}
h(f_k(x_k), f_k(y_k))\geqslant k \max\{\theta(|x_k-x_0|),
\theta(|y_k-x_0|)\}\,,\qquad k=1,2,\ldots\,.\end{equation}
Without loss of generality we may consider that
$|x_k-x_0|\geqslant|y_k-x_0|$ and, consequently, since $\theta$ is
an increasing function,
$$\max\{\theta(|x_k-x_0|),
\theta(|y_k-x_0|)\}=\theta(|x_k-x_0|)\,.$$
Thus, (\ref{eq2A}) implies that
\begin{equation}\label{eq2B}
h(f_k(x_k), f_k(y_k))\geqslant k \theta(|x_k-x_0|)\,,\qquad
k=1,2,\ldots\,.\end{equation}
By Lemma~\ref{lem2} there are paths $\gamma_k,$ $k=1,2,\ldots ,$
such that $\gamma_k$ lies in $\overline{B(x_0, r_k)}\cap D,$
$r_k=\max\{|x_k-x_0|, |y_k-x_0|\},$ $\gamma_k:[0, 1]\rightarrow
\overline{B(x_0, r_k)}\cap D,$ $\gamma_k(0)=x_k,$ $\gamma_k(1)=y_k,$
$\gamma_k(t)\in B(x_0, r_k)\cap D$ for $0<t<1,$ whenever
$|\gamma_k|:=\{x\in {\Bbb R}^n: \exists\,t\in [0, 1]:
\gamma_k(t)=x\}$ contains at most $m-1$ points in~$E.$

\medskip
Observe that, the path $f_k(\gamma_k)$ contains not more than $m-1$
points in $E_*.$ In the contrary case, there are at least $m$ such
points $b_{1}=f_k(\gamma_k(t_1)), b_{2}=f_k(\gamma_k(t_2)),\ldots,
b_m=f_k(\gamma_k(t_m)),$ $0\leqslant t_1\leqslant t_2\leqslant
\ldots\leqslant t_m\leqslant 1.$ But now the points
$a_{1}=\gamma_k(t_1), a_{2}=\gamma_k(t_2),\ldots, a_m=\gamma_k(t_m)$
are in $E=f_k^{\,-1}(E_*)$ and simultaneously belong to $\gamma_k.$
This contradicts the definition of $\gamma_k.$

\medskip
Let
\begin{gather*}b_{1}=f_k(\gamma_k(t_1)),\,
b_{2}=f_k(\gamma_k(t_2))\quad,\ldots,\quad
b_l=f_k(\gamma_k(t_l))\,,\\
0:=t_0\leqslant t_1\leqslant t_2\leqslant \ldots\leqslant
t_l\leqslant 1:= t_{l+1},\qquad 1\leqslant l\leqslant
m-1\,,\end{gather*}
be points in $f_k(\gamma_k)\cap E_*.$ By the triangle inequality,
\begin{equation}\label{eq6A}
h(f_k(x_k), f_k(y_k))\leqslant (m-1)\max_{0\leqslant r\leqslant
l}h(f_k(\gamma_k(t_{r})), f_k(\gamma_k(t_{r+1})))\,.
\end{equation}
Let $\max\limits_{0\leqslant r\leqslant
l}h(f_k(\gamma_k(t_{r})),f_k(\gamma_k(t_{r+1})))=
h(f_k(\gamma_k(t_{r(k)})),f_k(\gamma_k(t_{r(k)+1}))).$
Observe that, the set $G_k:=h(f_k(\gamma_k)|_{(t_{r(k)},
t_{r(k)+1})})$ belongs to $D_{f_k}^{\,\prime}\setminus E_*,$ because
it does not contain any point in $E_*,$ $D_{f_k}^{\,\prime}=f_k(D).$

\medskip
By Lemma~\ref{lem1}, a), $D\setminus E$ consists of finite number of
components $D_1,\ldots, D_s,$ $1\leqslant s<\infty.$ Thus, there is
a component $D_i,$ $1\leqslant i\leqslant s,$ such that a locus of
path $A_k:=|\gamma_k|_{(t_{r(k)}, t_{r(k)+1})}|$ belongs to $D_i$
for infinitely many $k.$ Without loss of generality we may assume
that the mentioned above holds for any $k\in {\Bbb N}$ and, in
addition, $A_k\subset D_1.$ By Lemma~\ref{lem1}, b) and c), $f_k$ is
an open, discrete and closed mapping of $D_1$ onto $K_k,$ where
$K_k$ is some a component of $D^{\,\prime}_{f_k}\setminus E_*.$

\medskip
By~(\ref{eq2B}) and (\ref{eq6A})
$$h(f_k(\gamma_k(t_{r(k)})), f_k(\gamma_k(t_{r(k)+1})))\geqslant
\frac{k}{m-1} \theta(|x_k-x_0|)$$
\begin{equation}\label{eq1A_1}
\geqslant \frac{k}{m-1} \max\{\theta(|\gamma_k(t_{r(k)})-x_0|),
\theta(|\gamma_k(t_{r(k)+1})-x_0|)\,.
\end{equation}
Here we have used that $\gamma_k$ lies in $B(x_0, r_k),$
$r_k=\max\{|x_k-x_0|, |y_k-x_0|\}=|x_k-x_0|,$ and that the function
$\theta(\varepsilon)$ is increasing by $\varepsilon.$ We set
$u_k:=\gamma_k(t_{r(k)}),$ $v_k:=\gamma_k(t_{r(k)+1}).$ Without loss
of generality we may assume that $|u_k-x_0|\geqslant |v_k-x_0|.$
Now, (\ref{eq1A_1}) implies that
\begin{equation}\label{eq1A_2}
h(f_k(u_k),f_k(v_k))\geqslant
\frac{k}{m-1}\theta(|u_k-x_0|)\,,\qquad k=1,2,\ldots \,.
\end{equation}
Fix $k\in {\Bbb N}$ and consider the function
$$\beta(x, y)=\frac{h(f_k(x),f_k(y))(m-1)}{\theta(|x-x_0|)}\,.$$
Observe that, the function $\beta$ is continuous at the point $(x,
y)=(u_k, v_k).$ Let $z_s$ and $w_s\in D_1\cap |\gamma_k|,$
$s=1,2,\ldots, $ be sequences converging to $u_k$ and $v_k,$
correspondingly. Now, given $\varepsilon=1$ there is a number
$s=s(k)$ such that $|z_s-u_k|<1/k,$ $|w_s-u_k|<1/k$ for $s\geqslant
s(k)$ and
$$|\beta(z_{s(k)}, w_{s(k)})-\beta(u_k, v_k)|<1\,.$$
From the latter relation, by the triangle inequality and
by~(\ref{eq1A_2}) we have that
$$\beta(z_{s(k)}, w_{s(k)})\geqslant \beta(u_k, v_k)-1\geqslant k-1\,.$$
Let $z_k:=z_{s(k)}$ and $w_k:=w_{s(k)}.$ Now, by the latter
inequality
\begin{equation}\label{eq1A_3}
h(f_k(z_k),f_k(w_k))\geqslant
\frac{k-1}{m-1}\theta(|z_k-x_0|)\,,\qquad k=1,2,\ldots \,.
\end{equation}
{\bf II.} Let us show that relation~(\ref{eq1A_3}) contradicts the
definition of class $\frak{R}^{E_*, E, N}_{\varphi, Q, \delta,
A_0}(D)$ for $n\geqslant 3$ (class $\frak{R}^{E_*, E, N}_{Q, \delta,
A_0}(D)$ for $n=2$) together with
relations~(\ref{eq3AB})--(\ref{eq3.7A}),
cf.~~\cite[Theorem~1.3]{Sev$_4$}.

\medskip
Let $K_{f_k}\subset K_k$ be a path connected continuum such that
$h(K_{f_k})\geqslant \delta$ and, in addition,
$h(f_k^{\,-1}(K_{f_k}),
\partial D)\geqslant \delta>0$ (such a continuum exists by the definition of
the class $\frak{R}^{E_*, E, N}_{\varphi, Q, \delta, A_0}(D)$
whenever $n\geqslant 3$ and the class $\frak{R}^{E_*, E, N}_{Q,
\delta, A_0}(D)$ whenever $n=2$). Now, $h((f_k^{\,-1}(K_{f_k}))\cap
D_1,
\partial D_1)\geqslant \delta>0.$ By
the definition, there is $\varepsilon_0>0$ such that
\begin{equation*}\label{eq18B}
(f_k^{\,-1}(K_{f_k})\cap D_1)\subset D_1\setminus B(x_0,
\varepsilon_0)\,.
\end{equation*}
Since $A_k:=|\gamma_k|_{(t_{r(k)}, t_{r(k)+1})}|\subset D_1\setminus
E$ and simultaneously $|\gamma_k|\subset B(x_0, r_k),$ the set $A_k$
belongs to some a component $Q_k$ of $(D\cap B(x_0, r_k))\setminus
E.$ Recall that, $z_k, w_k\in |\gamma_k|\cap D_1,$ consequently,
$z_k, w_k\in Q_k.$ By the condition 2b) we may join $z_k$ and
$w_k\in A_k$ in $Q_k$ by a path $\alpha_k$ in $B(x_0, r^{\,*}_k),$
where $r^{\,*}_k:=\max\{|z_k-x_0|, |w_k-x_0|\}.$ We may consider
that $|z_k-x_0|\geqslant |w_k-x_0|,$ so that $r^{\,*}_k=|z_k-x_0|.$
Observe that $|\alpha_k|\subset D\setminus E$ and simultaneously
$|\alpha_k|\cap D_1\ne\varnothing.$ Now, $|\alpha_k|\subset D_1\cap
B(x_0, r^{\,*}_k).$ Let $z_k^1, w_k^1\in K_{f_k}\subset
D^{\,\prime}_{f_k}$ be such that
\begin{equation}\label{eq8A}
h(K_{f_k})=h(z^1_k, w^1_k)\geqslant \delta\,.
\end{equation}
Since $K_{f_k}$ is path connected, we may join points $z^1_k, w^1_k$
by a Jordan path~$K^{\,\prime}_k$ inside $K_{f_k}.$ Due
to~(\ref{eq8A}) we obtain that
\begin{equation}\label{eq9A}
h(K^{\,\prime}_k)\geqslant h(z^1_k, w^1_k)\geqslant \delta\,.
\end{equation}
We also may consider that $K^{\,\prime}_k$ is Jordan. If the path
$f_k(|\alpha_k|))$ is not Jordan, we discard from $f_k(|\alpha_k|))$
no more than a countable number of its loops. Let $|L_k|\subset
f_k(|\alpha_k|))$ be a locus of the Jordan path $L_k,$ which is
obtained by such a rejection. By~\cite[Corollary~1.5.IV]{HW}, $L_k$
and $K^{\,\prime}_k$ do not split ${\Bbb R}^n$ for $n\geqslant 3,$
because, in this case, the set $L_k\cup |K^{\,\prime}_k|$ has a
topological dimension~1 as the union of two closed sets of
topological dimension~1 (see~\cite[Theorem~III 2.3]{HW}). Thus,
$R_k=R(|L_k|, |K^{\,\prime}_k|)$ is a ring in ${\Bbb R}^n.$

\medskip Now, let $n=2.$ We will reason as in proving Theorem~1.2 in~\cite{Sev$_4$}.
According to Antoine's theorem on the absence of wild arcs
(see~\cite[Theorem~II.4.3]{Keld}), there exists a homeomorphism
$\zeta:{\Bbb R}^2\rightarrow {\Bbb R}^2,$ which maps $|L_k|$ onto
some segment $I.$ It follows that, any points $x,y\in {\Bbb
R}^2\setminus |L_k|$ may be joined by a path $\gamma_k$ in ${\Bbb
R}^2\setminus |L_k|.$

Let us show that, any points $x,y\in {\Bbb R}^2\setminus
(|L_k|)\bigcup |K^{\,\prime}_k|)$ may be joined by a path $\gamma$
in ${\Bbb R}^2\setminus (|L_k|\bigcup |K^{\,\prime}_k|),$ as well.
Indeed, by the proving above, $x$ and $y$ may be joined by a path
$\gamma$ in ${\Bbb R}^2\setminus |L_k|).$ If $|\gamma|\cap
|K^{\,\prime}_k|=\varnothing,$ there is nothing to prove. In the
contrary case, let us consider a homeomorphism $\zeta:{\Bbb
R}^2\rightarrow {\Bbb R}^2,$ which maps $|K^{\,\prime}_k|$ onto some
segment $I.$ Let $\Pi$ be an open rectangular two of edges of which
are parallel to $I,$ and two of which are perpendicular to $I,$
while $I\subset \Pi.$ By the continuity of $\zeta,$ the sets
$\zeta(|L_k|)$ and $I=\zeta(|K_k^{\,\prime}|))$ are disjoint
compacta in ${\Bbb R}^2.$ Reducing $\Pi,$ we also may assume that
$\zeta(x)\not\in \Pi$ and $\zeta(y)\not\in \Pi.$ Thus, we may assume
that
\begin{equation}\label{eq1I}
\zeta(|L_k|)\cap\overline{\Pi}=\varnothing\,,\quad \zeta(x)\not\in
\Pi\,,\quad \zeta(y)\not\in \Pi\,.
\end{equation}
Let $\gamma:[0, 1]\rightarrow {\Bbb R}^2,$ $\gamma(0)=x,$
$\gamma(1)=y.$ Set
$$t_1:=\inf\limits_{t\in [0, 1], \zeta(\gamma(t))\in \Pi}t\,,\qquad t_2:=
\sup\limits_{t\in [0, 1], \zeta(\gamma(t))\in \Pi}t\,.$$
Since by the assumption $|\gamma|\cap |K^{\,\prime}_k|\ne
\varnothing,$ by~(\ref{eq1I}) and due
to~\cite[Theorem~1.I.5.46]{Ku$_2$}, we obtain that
$\zeta(\gamma(t_1))\in \partial \Pi$ and $\zeta(\gamma(t_2))\in
\partial \Pi.$ Now, we may replace a path $\gamma|_{[t_1, t_2]}$ by
a path $\alpha:[t_1, t_2]\rightarrow {\Bbb R}^2$ which does not
intersect $I.$ Finally, set
$$\widetilde{\gamma}(t)=\begin{cases}\gamma(t)\,, &t\in[0, 1]\setminus[t_1, t_2]\,, \\
\zeta^{\,-1}(\alpha(t))\,,& t\in [t_1, t_2]\end{cases}\,.$$
By the construction, $\widetilde{\gamma}$ does not intersect
$|L_k|,$ because $\gamma$ does not intersect $\zeta(|L_k|)$ by the
construction, in addition, $\zeta^{\,-1}(|\alpha|)\subset
\zeta^{\,-1}(\overline{\Pi})$ while
$\zeta^{\,-1}(\overline{\Pi})\cap |L_k|=\varnothing$ by the first
relation in~(\ref{eq1I}). On the other hand, $\widetilde{\gamma}$
does not intersect $|K^{\,\prime}_k|,$ because $\gamma(t)$ lies
outside some a neighborhood $U$ of $|K^{\,\prime}_k|$ for $t\in[0,
1]\setminus[t_1, t_2],$ and $\zeta^{\,-1}(\alpha(t))$ does not
intersect $\gamma$ by the construction of $\alpha.$ So,
$\widetilde{\gamma}$ is a required path. Therefore, $R=R(|L_k|,
|K^{\,\prime}_k|)$ is a ring domain both when $n$ is greater than or
equal to three and when $n=2.$

\medskip
In this case, we denote $\Gamma_k=\Gamma(|L_k|, |K^{\,\prime}_k|,
{\Bbb R}^n).$ Now, by Proposition~\ref{pr3_1}
\begin{equation}\label{eq10A}
M(\Gamma(|L_k|, |K^{\,\prime}_k|, {\Bbb
R}^n))\geqslant\frac{\omega_{n-1}}{\left(
\log\frac{2\lambda^2_n}{h(|K^{\,\prime}_k|)h(|L_k|)}\right)^{n-1}}\,.
\end{equation}
Due to~(\ref{eq9A}) and~(\ref{eq10A}), and by the definition of the
path~$K^{\,\prime}_k,$ we obtain that
\begin{equation}\label{eq15A}
M(\Gamma(|L_k|, |K^{\,\prime}_k|, {\Bbb
R}^n))\geqslant\frac{\omega_{n-1}}{\left(
\log\frac{2\lambda^2_n}{\delta\cdot h(f_k(z_k),
f_k(w_k))}\right)^{n-1}}\,.
\end{equation}
Since~$K_k$ is a $QED$-domain with a constant $A_0$
in~(\ref{eq4***}), by~(\ref{eq15A}) we obtain that
\begin{equation}\label{eq16A}
M(\Gamma(|L_k|, |K^{\,\prime}_k|,
K_k))\geqslant\frac{\omega_{n-1}}{A_0\left(
\log\frac{2\lambda^2_n}{\delta\cdot h(f_k(z_k),
f_k(w_k))}\right)^{n-1}}\,.
\end{equation}
Let $\Gamma^{\,*}_k$ be a family $\gamma:[0, 1)\rightarrow D_1$ of
all maximal $f_k$-liftings of paths $\gamma^{\,\prime}:[0,
1]\rightarrow K_k$ of the family $\Gamma_k=\Gamma(|L_k|,
|K^{\,\prime}_k|, K_k)$ starting at $|\alpha_k|.$ Such liftings
exist by Proposition~\ref{pr3}. By the same Proposition, due to the
closeness of~$f_k,$ we obtain that each path $\gamma\in
\Gamma^{\,*}_k$ has an extension $\gamma:[0, 1]\rightarrow D_1$ to
the point $b=1.$ Then $\gamma(1)\in f_k^{\,-1}(K_k)\subset
f_k^{\,-1}(K_{f_k}),$ that is,
\begin{equation}\label{eq4A}
\Gamma^{\,*}_k\subset \Gamma(|\alpha_k|, f_k^{\,-1}(K_{f_k})\cap
D_1, D_1)\,.
\end{equation}
Observe that,
\begin{equation}\label{eq3}
f_k(\Gamma^{\,*}_k)=\Gamma_k=\Gamma(|L_k|, |K^{\,\prime}_k|, K_k)\,.
\end{equation}
In this case, by the minorization of the modulus (see., e.g.,
\cite[Theorem~6.4]{Va}), by~(\ref{eq4A}) and~(\ref{eq3}) we obtain
that
\begin{equation}\label{eq20A}
M(f_k(\Gamma^{\,*}_k))=M(\Gamma(|L_k|, |K^{\,\prime}_k|, K_k))\,.
\end{equation}
Recall that $|\alpha_k|\in \overline{B(x_0, |z_k-x_0|)}$ and
$f_k^{\,-1}(K_{f_k})\cap D_1\subset D_1\setminus B(x_0,
\varepsilon_0).$ Now, by Lemmas~\ref{thOS4.1} and \ref{lem4} and
by~(\ref{eq20A}) for $n\geqslant 3$ (Lemmas~\ref{thOS4.1} and
\ref{lem3A} and by~(\ref{eq20A}) for $n=2$) we obtain that
\begin{equation}\label{eq21B}
M(\Gamma(|L_k|, |K^{\,\prime}_k|, K_k))\leqslant
N\int\limits_{A(x_0, |z_k-x_0|, \varepsilon_0)} Q(x)\cdot
\eta^n(|x-x_0|)\ dm(x)\,,
\end{equation}
where $\eta$ is an arbitrary nonnegative Lebesgue measurable
function satisfying the relation~(\ref{eq6B}) for $r_1=|z_k-x_0|,$
$r_2=\varepsilon_0.$ Put
$$\eta(t):=\begin{cases}\psi(t)/I(|z_k-x_0|, \varepsilon_0)\,,
& t\in(|z_k-x_0|, \varepsilon_0)\,,
\\
0\,,& t\not\in(|z_k-x_0|, \varepsilon_0)\,.
\end{cases}$$
Observe that, the function $\eta$ satisfies the
condition~(\ref{eq6B}) for $r_1=|z_k-x_0|,$ $r_2=\varepsilon_0.$
Now, by the relations~(\ref{eq21B}) and (\ref{eq3.7A}) we obtain
that
\begin{gather} M(\Gamma(|L_k|, |K^{\,\prime}_k|, K_k))\leqslant
\frac{N}{I^n(|z_k-x_0|, \varepsilon_0)}\int\limits_{A(x_0,
|z_k-x_0|, \varepsilon_0)} Q(x)\psi^n(|x-x_0|) \,dm(x)\nonumber\\
\label{eq21A} \leqslant K_0N\cdot I^{p-n}(|z_k-x_0|,
\varepsilon_0)\,.
\end{gather}
Combining~(\ref{eq16A}) and~(\ref{eq21A}), we obtain that
\begin{equation}\label{eq23A}
\frac{\omega_{n-1}}{A_0\left( \log\frac{2\lambda^2_n}{\delta\cdot
h(f_k(z_k), f_k(w_k))}\right)^{n-1}}\leqslant K_0N\cdot
I^{p-n}(|z_k-x_0|, \varepsilon_0)\,.
\end{equation}
By~(\ref{eq23A}), we obtain that
\begin{gather}
\log\frac{2\lambda^2_n}{\delta\cdot h(f_k(z_k), f_k(w_k))}\geqslant
(A_0K_0N\omega^{\,-1}_{n-1})^{-\frac{1}{n-1}}\cdot
I^{\frac{n-p}{n-1}}(|z_k-x_0|, \varepsilon_0)\,, \nonumber\\
\frac{2\lambda^2_n}{\delta\cdot h(f_k(z_k), f_k(w_k))}\geqslant
\exp\{(A_0K_0N\omega^{\,-1}_{n-1})^{-\frac{1}{n-1}}\cdot
I^{\frac{n-p}{n-1}}(|z_k-x_0|, \varepsilon_0)\}\,, \nonumber\\
h(f_k(z_k), f_k(w_k))\leqslant
\frac{2\lambda^2_n}{\delta}\exp\{-(A_0K_0N\omega^{\,-1}_{n-1})^{-\frac{1}{n-1}}\cdot
I^{\frac{n-p}{n-1}}(|z_k-x_0|, \varepsilon_0)\}\nonumber \\
\label{eq25A} =\frac{2\lambda^2_n}{\delta}\theta(|z_k-x_0|)\,.
\end{gather}
It follows by~(\ref{eq25A}) that the relation $$\frac{h(f_k(z_k),
f_k(w_k))}{\theta(|z_k-x_0|)}$$ is upper bounded for sufficiently
large $k\in {\Bbb N}.$ The latter contradicts with~(\ref{eq1A_3}).
It finishes the proof.~$\Box$ \end{proof}

\medskip
We also present an analogue of Lemma~\ref{lem3} for domains with a
more complex boundary structure, for which it is necessary to
involve prime ends.

\medskip
\begin{lemma}\label{lem5} {\it\, Assume
that, under conditions of Theorem~\ref{th5}, instead of the
relation~(\ref{eq1EA}) the following is true: for every
$x_0\in\partial D$ there are $\varepsilon_0>0,$ $1\leqslant p<n,$
$K_0>0,$ and a Lebesgue measurable function $\psi:(0,
\varepsilon_0)\rightarrow [0, \infty]$ such that the
relations~(\ref{eq3AB})--(\ref{eq3.7A}) hold.
Then for each $P\in E_D$ there exists $y_0\in \partial D$ such that
$I(P)=y_0$ while $I(P)$ denotes the impression of $P,$
see~(\ref{eq1A_7}). In addition, there exists a neighborhood $U$ of
$P_0$ in $(\overline{D}_P, \rho)$ and $\alpha>0$ such that the
inequality
\begin{gather}\label{eq2A_2}
h(f(x), f(y))\leqslant \alpha\cdot\max\{\theta(|x-x_0|),
\theta(|y-x_0|)\}\,,\\
\nonumber
\theta(\varepsilon):=\exp\{-(A_0K_0N\omega^{\,-1}_{n-1})^{-\frac{1}{n-1}}\cdot
I^{\frac{n-p}{n-1}}(\varepsilon, \varepsilon_0)\}\,,
\end{gather}
holds for any $x, y\in U\cap D$ and all $f\in\frak{R}^{E_*, E,
N}_{\varphi, Q, \delta, A_0}(D),$ $x_0:=I(P_0),$ whenever
$n\geqslant 3$ ($f\in\frak{R}^{E_*, E, N}_{\varphi, Q, \delta,
A_0}(D)$ for $n=2$). }
\end{lemma}

\medskip
\begin{proof} Observe that, $E$ is closed in $D$ due to the
relation $f^{\,-1}(E_*)=E,$ because while $E_*$ is closed in
$\overline{{\Bbb R}^n}$ by the assumption, see e.g.
\cite[Theorem~1.IV.13.1]{Ku$_1$}. Since the set $B(y_0, r)\cap D$ is
finitely connected for all $0<r<r_0,$ the domain $D$ is finitely
connected on its boundary. Therefore, the domain $D$ is uniform in
the sense of N\"akki (see \cite[Theorem~3.2]{Na$_1$}). In other
words, for every $r>0$ there exists a constant $\delta>0$ such that
the inequality  $M(\Gamma(F^{\,*},F, D))\geqslant \delta$ holds for
all continua $F, F^*\subset D$ such that $h(F)\geqslant r$ and
$h(F^{\,*})\geqslant r.$ Arguing similarly to the proof of
Theorem~1.3 in \cite{Sev$_4$}, we may show that for any $P\in E_D$
there exists $y_0\in
\partial D$ such that $I(P)=y_0.$ Now we will proceed directly to the proof of
relation~(\ref{eq2A_2}). In general, we will adhere to the scheme of
proof of Lemma~\ref{lem3}, cf. proof of Theorem~1.3 in
\cite{Sev$_4$}.

\medskip
As under the proof of Lemma~\ref{lem1}, it is sufficiently to prove
Lemma~\ref{lem5} for $x, y\in D\setminus E,$ because the set $E$ is
nowhere dense in $D.$ Let us prove by contradiction, i.e., assume
that the statement of Theorem~\ref{th5} is not true for $x, y\in
D\setminus E.$ Now, there are sequences $x_k, y_k\rightarrow P_0$ as
$k\rightarrow\infty,$ $x_k, y_k\in D\setminus E,$ and $f_k\in
\frak{R}^{E_*, E, N}_{\varphi, Q, \delta, A_0}(D)$ whenever
$n\geqslant 3$ ($f_k\in\frak{R}^{E_*, E, N}_{Q, \delta, A_0}(D)$ for
$n=2$) such that
\begin{equation}\label{eq2A_1}
h(f_k(x_k), f_k(y_k))\geqslant k \max\{\theta(|x_k-x_0|),
\theta(|y_k-x_0|)\}\,,\qquad k=1,2,\ldots\,,\end{equation}
where $x_0=I(P_0).$ Without loss of generality we may consider that
$|x_k-x_0|\geqslant |y_k-x_0|$ and, consequently,
$\theta(|x_k-x_0|)\geqslant \theta(|y_k-x_0|)$ for any $k\in {\Bbb
N}.$ Thus, (\ref{eq2A_1}) implies that
\begin{equation}\label{eq2B_1}
h(f_k(x_k), f_k(y_k))\geqslant k \theta(|x_k-x_0|)\,,\qquad
k=1,2,\ldots\,.\end{equation}
We may consider that, $x_k, y_k\in d_k,$ $k=1,2,\ldots, $ where
$d_k$ is a sequence of decreasing domains in $P_0.$ Let
$\varepsilon_0>0$ be such that $h(f_k^{\,-1}(K_{f_k}),
\partial D)\geqslant \varepsilon_0>0$ for any $k\in {\Bbb N}$ (such $\varepsilon_0>0$ exists
by the assumption of the lemma). Since $I(P_0)=x_0,$ there exists
$k_0\in {\Bbb N}$ such that $d_k\subset B(x_0, \varepsilon_0)$ for
any $k\geqslant k_0.$ Now, let $P_k$ be a component of $B(x_0,
\varepsilon_0)\cap D$ containing $d_k.$ Since $|x_k-x_0|\geqslant
|y_k-x_0|$ for any $k\in {\Bbb N}$ by the assumption, $x_k, y_k\in
\overline{B(x_0, |x_k-x_0|)}.$ Now, by the assumption 2a) we may
join $x_k$ and $y_k$ by a path $\gamma_k$ in $P_k\cap
\overline{B(x_0, |x_k-x_0|)},$ $\gamma_k(0)=x_k,$ $\gamma_k(1)=y_k,$
$\gamma_k(t)\in P_k\cap B(x_0, r_k)$ for $0<t<1,$ $r_k:=|x_k-x_0|,$
whenever $|\gamma_k|:=\{x\in {\Bbb R}^n: \exists\,t\in [0, 1]:
\gamma_k(t)=x\}$ contains at most $m-1$ points in~$E.$

\medskip
Similar to what was done in the proof of Lemma~\ref{lem3}, we may
show that the path $f_k(\gamma_k)$ contains not more than $m-1$
points in $E_*.$ Let
\begin{gather*}b_{1}=f_k(\gamma_k(t_1)),\,
b_{2}=f_k(\gamma_k(t_2))\quad,\ldots,\quad
b_l=f_k(\gamma_k(t_l))\,,\\
0:=t_0\leqslant t_1\leqslant t_2\leqslant \ldots\leqslant
t_l\leqslant 1:= t_{l+1},\qquad 1\leqslant l\leqslant
m-1\,,\end{gather*}
be points in $f_k(\gamma_k)\cap E_*.$ By the triangle inequality,
\begin{equation}\label{eq6A_1}
h(f_k(x_k), f_k(y_k))\leqslant (m-1)\max_{0\leqslant r\leqslant
l}h(f_k(\gamma_k(t_{r})), f_k(\gamma_k(t_{r+1})))\,.
\end{equation}
Let $\max\limits_{0\leqslant r\leqslant
l}h(f_k(\gamma_k(t_{r})),f_k(\gamma_k(t_{r+1})))=
h(f_k(\gamma_k(t_{r(k)})),f_k(\gamma_k(t_{r(k)+1}))).$
Observe that, the set $G_k:=h(f_k(\gamma_k)|_{(t_{r(k)},
t_{r(k)+1})})$ belongs to $D_{f_k}^{\,\prime}\setminus E_*,$ because
it does not contain any point in $E_*.$

\medskip
Observe that, $(B(x_0, r)\cap D)\setminus E$ consists at most a
finite number of components for sufficiently small $r>0.$ Indeed,
$B(x_0, r)\cap D$ consists at most of $p$ components $K_1, K_2,
\ldots, K_j,$ $1\leqslant j\leqslant p,$ by the condition~2a). By
the condition~2b) $(B(x_0, r)\cap K_i)\setminus E$ consists at most
of $m$ components $K_{i1},\ldots, K_{im_i},$ $1\leqslant
m_i\leqslant m,$ for any $1\leqslant i\leqslant p.$ Now, $(B(x_0,
r)\cap D)\setminus
E=\bigcup\limits_{i=1}^j\bigcup\limits_{l=1}^{m_i}K_{il}.$ By
Lemma~\ref{lem1}, a), $D\setminus E$ consists of finite number of
components $D_1,\ldots, D_s,$ $1\leqslant s<\infty.$ Thus, there is
a component $D_i,$ $1\leqslant i\leqslant s,$ such that a locus of
path $A_k:=|\gamma_k|_{(t_{r(k)}, t_{r(k)+1})}|$ belongs to $D_i$
for infinitely many $k.$ Without loss of generality we may assume
that the mentioned above holds for any $k\in {\Bbb N}$ and, in
addition, $A_k\subset D_1.$ Lemma~\ref{lem1}, b) and c), $f_k$ is an
open, discrete and closed mapping of $D_1$ onto $K_k,$ where $K_k$
is some a component of $D^{\,\prime}_{f_k}\setminus E_*.$

\medskip
By~(\ref{eq2B_1}) and (\ref{eq6A_1})
\begin{gather}
h(f_k(\gamma_k(t_{r(k)})),f_k(\gamma_k(t_{r(k)+1})))\nonumber\\
\label{eq1A_6}\geqslant \frac{k}{m-1} k \theta(|x_k-x_0|)\geqslant
\frac{k}{m-1} \max\{\theta(|\gamma_k(t_{r(k)})-x_0|)|,
\theta(|\gamma_k(t_{r(k)+1})-x_0|)|\}\,.
\end{gather}
Here we have used that $\gamma_k$ lies in $B(x_0, r_k),$
$r_k=|x_k-x_0|.$ We set $u_k:=\gamma_k(t_{r(k)}),$
$v_k:=\gamma_k(t_{r(k)+1}).$ Without loss of generality we may
assume that $|u_k-x_0|\geqslant |v_k-x_0|.$ Now, (\ref{eq1A_6})
implies that
\begin{equation}\label{eq1A_5}
h(f_k(u_k),f_k(v_k))\geqslant \frac{k}{m-1}
\theta(|u_k-x_0|)\,,\qquad k=1,2,\ldots \,.
\end{equation}
Arguing similarly to the proof of Lemma~\ref{lem3}, we may conclude
that~(\ref{eq1A_5}) implies the existing of the sequences $z_k,
w_k\in D_1\cap |\gamma_k|,$ $z_k, w_k\rightarrow x_0$ as
$k\rightarrow\infty$ such that
\begin{equation}\label{eq1A_4}
h(f_k(z_k),f_k(w_k))\geqslant
\frac{k-1}{m-1}\theta(|z_k-x_0|)\,,\qquad k=1,2,\ldots \,.
\end{equation}
Let us show that relation~(\ref{eq1A_4}) contradicts the definition
of class $\frak{R}^{E_*, E, N}_{\varphi, Q, \delta, A_0}(D)$
whenever $n\geqslant 3$ (class $\frak{R}^{E_*, E, N}_{Q, \delta,
A_0}(D)$ for $n=2$) together with relation~(\ref{eq1EA}).

\medskip
Indeed, by the definition of $\varepsilon_0>0$ mentioned above,
\begin{equation*}\label{eq18A}
(f_k^{\,-1}(K_{f_k})\cap D_1)\subset D_1\setminus B(x_0,
\varepsilon_0)\,.
\end{equation*}
Since $A_k:=|\gamma_k|_{(t_{r(k)}, t_{r(k)+1})}|\subset D_1$ and
simultaneously $|\gamma_k|\subset P_k\cap B(x_0, r_k)\subset B(x_0,
r_k),$ the set $A_k$ belongs one and only one component $Q_k$ of
$D\cap B(x_0, r_k).$ Let $\widetilde{Q}_k$ be a component of
$Q_k\setminus E$ containing $A_k.$ (Such a component
$\widetilde{Q}_k$ exists because $A_k$ is connected and $A_k\subset
Q_k\setminus E$). Now, by the assumption~2b), we may join $z_k$ and
$w_k$ in $\widetilde{Q}_k\cap \overline{B(x_0, r^{\,*}_k)}$ by a
path $\alpha_k,$ where $r^{\,*}_k:=\max\{|z_k-x_0|, |w_k-x_0|\}.$ We
may consider that $|z_k-x_0|\geqslant |w_k-x_0|,$ so that
$r^{\,*}_k=|z_k-x_0|.$ Observe that, $\alpha_k\subset
\widetilde{Q}_k\subset D_1.$ Indeed, $\alpha_k\subset
\widetilde{Q}_k$ by the construction. In addition, $z_k, w_k\in D_1$
by the construction, so that $\widetilde{Q}_k\cap
D_1\ne\varnothing.$ Since the connected set $\widetilde{Q}_k$ in
$D\setminus E$ may belongs to one and only one component of
$D\setminus E,$ $\widetilde{Q}_k\subset D_1,$ as required.

\medskip
Let $z_k^1, w_k^1\in K_{f_k}\subset D^{\,\prime}_{f_k}$ be such that
\begin{equation}\label{eq8A_1}
h(K_{f_k})=h(z^1_k, w^1_k)\geqslant \delta\,.
\end{equation}
Since $K_{f_k}$ is path connected, we may join points $z^1_k, w^1_k$
by a Jordan path~$K^{\,\prime}_k$ inside $K_{f_k}.$ Due
to~(\ref{eq8A_1}) we obtain that
$$
h(K^{\,\prime}_k)\geqslant h(z^1_k, w^1_k)\geqslant \delta\,.
$$
We also may consider that $K^{\,\prime}_k$ is Jordan. If the path
$f_k(|\alpha_k|))$ is not Jordan, we discard from $f_k(|\alpha_k|))$
no more than a countable number of its loops. Let $|L_k|\subset
f_k(|\alpha_k|))$ be a locus of the Jordan path $L_k,$ which is
obtained by such a rejection. Observe that, $L_k$ and
$K^{\,\prime}_k$ do not split ${\Bbb R}^n.$ Indeed, for $n\geqslant
3$ the set $L_k\cup |K^{\,\prime}_k|$ has a topological dimension~1
as the union of two closed sets of topological dimension~1
(see~\cite[Theorem~III 2.3]{HW}). Then $L_k\cup |K^{\,\prime}_k|$
does not split ${\Bbb R}^n$ (see~\cite[Corollary~1.5.IV]{HW}). If
$n=2,$ we are reasoning similarly to the proof of Lemma~\ref{lem3}.
Now
$R_k=R(|L_k|, |K^{\,\prime}_k|)$ is a ring in ${\Bbb R}^n.$

In this case, we denote $\Gamma_k=\Gamma(|L_k|, |K^{\,\prime}_k|,
{\Bbb R}^n).$ Now, by Proposition~\ref{pr3_1}
\begin{equation*}\label{eq10A_1}
M(\Gamma(|L_k|, |K^{\,\prime}_k|, {\Bbb
R}^n))\geqslant\frac{\omega_{n-1}}{\left(
\log\frac{2\lambda^2_n}{h(|K^{\,\prime}_k|)h(|L_k|)}\right)^{n-1}}\,.
\end{equation*}
Reasoning further in the same way as after relation~(\ref{eq10A}) in
the proof of Lemma~\ref{lem3}, we arrive at the relation
\begin{gather*}h(f_k(z_k), f_k(w_k))\leqslant
\frac{2\lambda^2_n}{\delta}\exp\{-(A_0K_0N\omega^{\,-1}_{n-1})^{-\frac{1}{n-1}}\cdot
I^{\frac{n-p}{n-1}}(|z_k-x_0|, \varepsilon_0)\}\\
=\frac{2\lambda^2_n}{\delta}\theta(|z_k-x_0|)\,.
\end{gather*}
The latter contradicts with~(\ref{eq1A_4}). It finishes the
proof.~$\Box$ \end{proof}

\section{Proof of Theorems~\ref{th4} and \ref{th5}}

{\it Proof of Theorem~\ref{th4}.} In Lemma~\ref{lem3}, we set
$\psi(t)=\frac{1}{t}.$ Now, $I(\varepsilon,
\varepsilon_0)=\log\frac{\varepsilon_0}{\varepsilon}.$ In this case,
Proposition~\ref{pr1} implies that the relation~(\ref{eq3.7A}) holds
with some $K_0=C_1$ and $p=1.$ Let $|x-x_0|\geqslant |y-x_0|.$ Now,
by Lemma~\ref{lem3}
$$
h(f(x), f(y))\leqslant \alpha
\exp\left\{-(A_0K_0N)^{-\frac{1}{n-1}}\cdot
\log\frac{\varepsilon_0}{|x-x_0|}\right\}=\alpha\left(\frac{|x-x_0|}{\varepsilon_0}\right)
^{(A_0K_0N)^{-\frac{1}{n-1}}}\,.~\Box$$

\medskip
{\it Proof of Theorem~\ref{th5}} is similar to the proof of
Theorem~\ref{th4} and based on Lemma~\ref{lem5}.$~\Box$

\medskip
\begin{example}\label{ex1B}
Let $D={\Bbb D}\setminus I,$ let ${\Bbb D}=\{z\in \Bbb C: |z|<1\},$
let $I=\{z\in \Bbb C: z=t, t\in [0, 1)\}$ and let $f(z)=(f_2\circ
f_1)(z),$ where $f_1(z)=|z|r^{i\varphi/2}$ is the main branch of
$\sqrt{z}$ and $f_2(z)=z^3,$ see Figure~\ref{fig7A}.%
\begin{figure}
\centering\includegraphics[width=300pt]{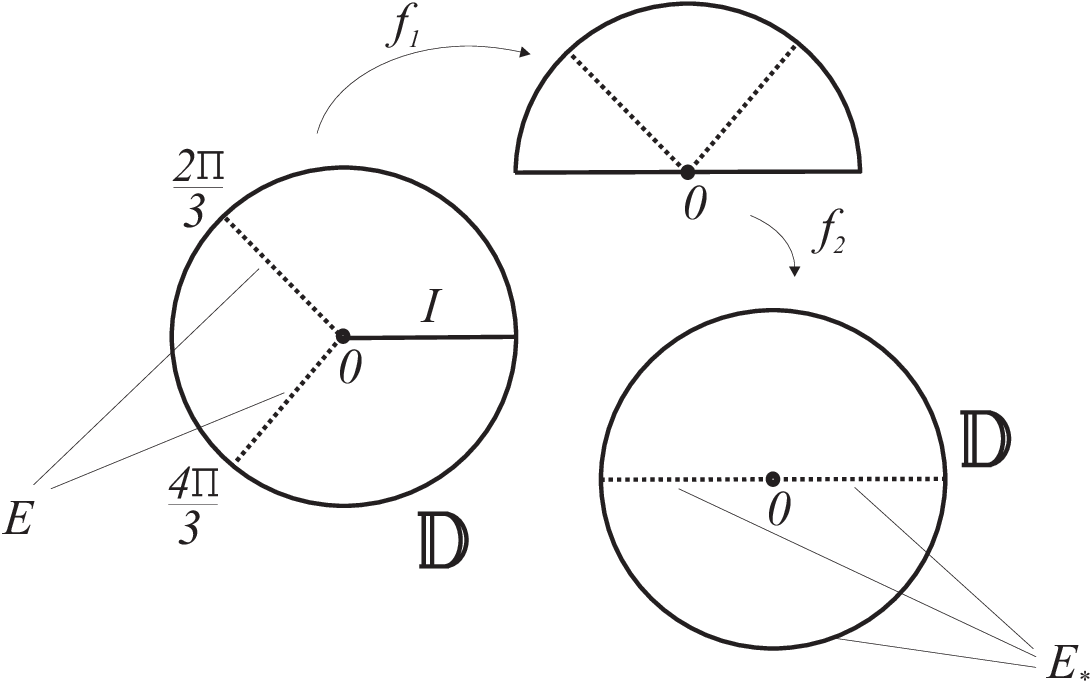}
\caption{Example~\ref{ex1B}.}\label{fig7A}
\end{figure}
\end{example}
We verify all the conditions of Theorem~\ref{th5} with $n=2.$
Indeed, $D$ is regular domain by the Riemannian mapping theorem. It
is clear that $f(D)={\Bbb D}$ and $E_*:=C(f, \partial D)={\Bbb
S}^1\cup I_1,$ where $I_1=\{z\in \Bbb C: z=t, t\in (-1, 1)\}$ and
${\Bbb S}^1=\{z\in {\Bbb C}: |z|=1\}.$ Besides that,
$E=f^{\,-1}(I_1)=J_1\cup J_2,$ where $J_1=\{z\in \Bbb C:
z=re^{i\frac{2\pi}{3}}, r\in (0, 1)\}$ and $J_2=\{z\in \Bbb C:
z=re^{i\frac{4\pi}{3}}, r\in (0, 1)\}.$ Obviously, $E$ is nowhere
dense in $D.$ Observe that $D$ is finitely connected on $E\cup
\partial D,$ i.e., for any $z_0\in E\cup\partial D$ and any
neighborhood $\widetilde{U}$ of $z_0$ there is a neighborhood
$\widetilde{V}\subset \widetilde{U}$ of $z_0$ such that $(D\cap
\widetilde{V})\setminus E$ consists of finite number of components.
Indeed, the above $V$ may be chosen such that $(D\cap
\widetilde{V})\setminus E$ consists of one component whenever
$z_0\in
\partial
D\setminus\{I\cup\{e^{i\frac{2\pi}{3}}\}\cup\{e^{i\frac{4\pi}{3}}\}\};$
of two components whenever $z_0\in E\cup \{e^{i\frac{2\pi}{3}}\}\cup
\{e^{i\frac{4\pi}{3}}\}\cup I\setminus\{0\};$ and of three
components whenever $z_0=0.$

\medskip
Obviously, all of the conditions in 2), Theorem~\ref{th5}, are
satisfied, because, firstly, $D$ forms at most two connected
components with intersection of infinitesimal boundary balls $B(y_0,
r),$ $y_0\in {\Bbb S}.$ Besides that, each component $K$ of $B(y_0,
r)\cap D$ is convex for sufficiently small $r>0.$ Therefore, the
corresponding condition ``$x, y\in K$ may be joined by a path
$\gamma:[a, b]\rightarrow {\Bbb R}^n$ such that $|\gamma|\in K\cap
\overline{B(y_0 , \max\{|x-y_0|, |y-y_0|\})}$'' obviously holds.
Secondly, $K\setminus E$ coincide with $K,$ excluding the points
$y_0\in \{0, e^{i\frac{2\pi}{3}}, e^{i\frac{4\pi}{3}}\}.$ In the
latter case, $K\setminus E$ consists at most three components (two
components for $y_0=e^{i\frac{2\pi}{3}}$ or
$y_0=e^{i\frac{4\pi}{3}},$ and three components for $y_0=0$). Each
of these components are convex for small $r>0.$ It means that the
above condition ``for each component $K_1$ of the set $K\setminus E$
the following condition is fulfilled: any $x, y\in K_1$ may be
joined by a path $\gamma:[a, b]\rightarrow {\Bbb R}^n$ such that
$|\gamma|\in K_1\cap \overline{B(y_0 , \max\{|x-y_0|, |y-y_0|\})}$''
holds, as well. The relation~(\ref{eq1EA}) (condition~3)) holds
because $K_I(x, f)=1.$

Finally, we show that $f\in\frak{R}^{E_*, E, N}_{Q, \delta, A_0}(D)$
for some parameters $E_*, E, N, Q, \delta$ and $A_0.$ Indeed, the
conditions 1) $C(f, \partial D)\subset E_*$ and 3) $f^{\,-1}(E_*)=E$
in the definition of the class $\frak{R}^{E_*, E, N}_{Q, \delta,
A_0}(D)$ are fulfilled by the construction of $E$ and $E_*$
mentioned above. The condition 2) ``for each component $K$ of
$D^{\,\prime}_f\setminus E_*,$ $D^{\,\prime}_f:=f(D),$ there is a
continuum $K_f\subset K$ such that $h(K_f)\geqslant \delta$ and
$h(f^{\,-1}(K_f), \partial D)\geqslant \delta>0$'' trivially holds
for any component $K$ of $D^{\,\prime}_f={\Bbb D}$ and any continuum
in $K.$ The condition 4) ``any component $G$ of $f(D)\setminus E_*$
satisfies the relation~(\ref{eq4***})'' holds because $f(D)\setminus
E_*={\Bbb D}\setminus E_*$ consists of two components
$\widetilde{D}_1={\Bbb D}_+:=\{z=x+iy\in {\Bbb D}: y>0\}$ and
$\widetilde{D}_2={\Bbb D}_-:=\{z=x+iy\in {\Bbb D}: y<0\},$ each of
which are $QED$-domains (see, e.g., \cite[Theorem~2.22]{GM}).
Obviously, $f$ is bounded, open and discrete mapping in $W_{\rm
loc}^{1, 1}(D)$ such that $K_I(x, f)\leqslant 1$ a.e. and $N(f,
D)\leqslant 3.$

\medskip
\begin{example}
Let
$$\beta(t)= \left\{
\begin{array}{rr}
1, & t\in\bigcup\limits_{k\geqslant 1}\left[\frac{k}{k+1}, \frac{k+1}{k+2}-2^{-4k-1}\right],\\
2^k,  &
t\in\left(\frac{k+1}{k+2}-2^{-4k-1},\,\frac{k+1}{k+2}\right)\,,\textrm{
for some}\,\,k>0
\end{array}
\right. $$ $k=0,1,2,\ldots .$
$Q(z)=\beta(|z|).$ Observe that, the function $Q$
satisfies~(\ref{eq1EA}) at any $z_0\in \partial D$ (see
\cite[Example~2.2]{MSS}). We can see that the mapping
$$w=g(re^{i\theta})=e^{i\theta+\int\limits_{1}^r(\beta(t)/t)\,dt}$$
is a homeomorphism of the Sobolev class $W_{\rm loc}^{1, 1}({\Bbb
C})$ and has a dilatation $K_I(z, g)$ equal to $\beta(|z|).$ Note
that this map takes the unit disk onto itself; moreover, $g|_{{\Bbb
S}^1}$ is H\"{o}lder continuous with an arbitrary exponent on ${\Bbb
S}^1$ since $f(z)=z$ at the points $z\in{\Bbb S}^1.$ Observe that,
$g$ maps $E$ onto $E,$ where $E$ is a set from Example~\ref{ex1B}.
Now, the mapping $F=f\circ g$ satisfies all the conditions of
Theorem~\ref{th5} with $D, E_*, E, N, \delta$ and $A_0$ mentioned in
this example, while $Q$ is defined above. In this example, it is
important that the function $Q$ is unbounded, although
condition~(\ref{eq1EA}) is still satisfied. Furthermore, the mapping
$F$ does not preserve the boundary of the domain, since the mapping
$g$ is a homeomorphism, and the mapping $f$ also does not preserve
the boundary.
\end{example}

\section{Some another distance distortion estimates}

Set
\begin{equation*}\label{eq12_2}
q_{y_0}(r)=\frac{1}{\omega_{n-1}r^{n-1}}\int\limits_{S(y_0,
r)}Q(y)\,d\mathcal{H}^{n-1}(y)\,, \end{equation*}
and $\omega_{n-1}$ denotes the area of the unit sphere ${\Bbb
S}^{n-1}$ in ${\Bbb R}^n.$ The following statement holds.

\medskip
\begin{theorem}\label{th1}
{\it\, Assume that, under conditions of Theorem~\ref{th4}, instead
of the relation~(\ref{eq1E_1}) the following is true:
\begin{equation}\label{eq13}
\int\limits_{\varepsilon}^{\delta_0}\frac{dt}
{tq_{x_0}^{\,\frac{1}{n-1}}(t)}<\infty\qquad\forall\quad
\varepsilon\in (0, \delta_0)\,,\qquad
\int\limits_{0}^{\delta_0}\frac{dt}
{tq_{x_0}^{\,\frac{1}{n-1}}(t)}=\infty
\end{equation}
for some $\delta_0=\delta_0(x_0)>0$ and sufficiently small
$0<\varepsilon<\delta_0.$ Then
$$
h(f(x), f(y))\leqslant
\alpha\cdot\exp\left\{-(A_0N)^{-\frac{1}{n-1}}
\int\limits_{|x-x_0|}^{\varepsilon_0}\frac{dt}
{tq_{x_0}^{\frac{1}{n-1}}(t)}\right\}$$
holds for any $x, y\in B(x_0, \widetilde{\varepsilon}(x_0))\cap D$
and all $f\in\frak{R}^{E_*, E, N}_{\varphi, Q, \delta, A_0}(D)$
whenever $n\geqslant 3$ ($f\in\frak{R}^{E_*, E, N}_{Q, \delta,
A_0}(D)$ for $n=2$).}
\end{theorem}

\medskip
\begin{proof}
We set $\varepsilon_0=\delta_0,$ where $\delta_0$ is a number
from~(\ref{eq13}). The relation~(\ref{eq13}) implies that
$I(\varepsilon, \varepsilon_0)>0$ for all $\varepsilon\in(0,
\varepsilon_1)$ and some $\varepsilon_1\in (0, \varepsilon_0).$ Let
$0<\varepsilon<\varepsilon_1.$ Let us consider a function
$I(\varepsilon, \varepsilon_0)=\int\limits
_{\varepsilon}^{\varepsilon_0}\psi(t)dt,$ where
%
%\begin{equation}\label{eq44**}
$$\psi(t)\quad=\quad \left \{\begin{array}{rr}
1/[tq^{\frac{1}{n-1}}_{x_0}(t)]\ , & \ t\in (\varepsilon,
\varepsilon_0)\ ,
\\ 0\ ,  &  \ t\notin (\varepsilon,
\varepsilon_0)\ ,
\end{array} \right.$$
%\end{equation}
%
(here we use relations $a/\infty = 0$ for $a\ne\infty,$ $a/0=\infty
$ for $a>0$ and $0\cdot\infty =0,$ see \cite[Ch.~I]{Sa}). Observe
that, $\psi$ satisfies the relation~(\ref{eq3.7A}). Indeed, by the
Fubini theorem (see e.g.~\cite[Theorem~8.1.III]{Sa}),
$$
\int\limits_{\varepsilon<|x-x_0|<\varepsilon_0}
Q(x)\cdot\psi^n(|x-x_0|)\,dm(x)=\omega_{n-1}\cdot
\int\limits_{\varepsilon}^{\varepsilon_0}\frac{dt}
{tq_{x_0}^{\frac{1}{n-1}}(t)}=\omega_{n-1}\cdot I(\varepsilon,
\varepsilon_0)\,.
$$
Thus the relation~(\ref{eq3.7A}) holds for $p=1$ and
$K_0=\omega_{n-1}.$ Let $|x-x_0|\geqslant |y-x_0,|$ Now, by
Lemma~\ref{lem3}
$$
h(f(x), f(y))\leqslant
\alpha\cdot\exp\left\{-(A_0N)^{-\frac{1}{n-1}}
\int\limits_{|x-x_0|}^{\varepsilon_0}\frac{dt}
{tq_{x_0}^{\frac{1}{n-1}}(t)}\right\}$$
holds for some $\alpha>0,$ any $x, y\in B(x_0,
\widetilde{\varepsilon}(x_0))\cap D$ and all $f\in\frak{R}^{E_*, E,
N}_{\varphi, Q, \delta, A_0}(D)$ whenever $n\geqslant 3$
($f\in\frak{R}^{E_*, E, N}_{Q, \delta, A_0}(D)$ for $n=2$), as
required.~$\Box$
\end{proof}

\begin{theorem}\label{th2}
{\it\, Assume that, under conditions of Theorem~\ref{th5}, instead
of the relation~(\ref{eq1EA}) the relation~(\ref{eq13}) is true for
any $x_0\in I(P_0),$ some $\delta_0=\delta_0(x_0)>0$ and
sufficiently small $0<\varepsilon<\delta_0.$ Then there is
$0<\widetilde{\varepsilon}_0<\varepsilon_0$ and $\alpha>0,$ does not
depending on $f,$ such that
$$
h(f(x), f(y))\leqslant
\alpha\cdot\exp\left\{-(A_0N)^{-\frac{1}{n-1}}
\int\limits_{|x-x_0|}^{\varepsilon_0}\frac{dt}
{tq_{x_0}^{\frac{1}{n-1}}(t)}\right\}$$
holds for some $\alpha>0,$ $x_0=I(P_0),$ some a neighborhood $U$ of
$P_0,$ any $x, y\in U\cap D$ and all $f\in\frak{R}^{E_*, E,
N}_{\varphi, Q, \delta, A_0}(D)$ whenever $n\geqslant 3$
($f\in\frak{R}^{E_*, E, N}_{Q, \delta, A_0}(D)$ for $n=2$).}
\end{theorem}

\medskip
{\it Proof of Theorem~\ref{th2}} is based on Lemma~\ref{lem5} and is
similar to the proof of Theorem~\ref{th1}.~$\Box$

\medskip
In accordance with~\cite{MRSY}, we say that a function
$\vartheta:D\rightarrow{\Bbb R}$ has a {\it finite mean oscillation}
at a point $x_0\in D,$ write $\vartheta\in FMO(x_0),$ if
$$\limsup\limits_{\varepsilon\rightarrow
0}\frac{1}{\Omega_n\varepsilon^n}\int\limits_{B( x_0,\,\varepsilon)}
|\vartheta(x)-\overline{\vartheta}_{\varepsilon}|\ dm(x)<\infty\,,
$$
where $\overline{\vartheta}_{\varepsilon}=\frac{1}
{\Omega_n\varepsilon^n}\int\limits_{B(x_0,\,\varepsilon)}
\vartheta(x) \,dm(x)$ and $\Omega_n$ is the volume of the unit ball
${\Bbb B}^n$ in ${\Bbb R}^n.$
We also say that a function $\vartheta:D\rightarrow{\Bbb R}$ has a
finite mean oscillation at $A\subset \overline{D},$ write
$\vartheta\in FMO(A),$ if $\vartheta$ has a finite mean oscillation
at any point $x_0\in A.$ The following statement holds, see e.g.
\cite[Corollary~6.3]{MRSY}.

\medskip
\begin{proposition}\label{pr2}
{\it\, If $Q:{\Bbb R}^n\rightarrow[0, \infty],$ $Q\in FMO(x_0),$
then
$$\int\limits_{A(x_0, \varepsilon, \varepsilon_0)}\frac{Q(x)\,dm(x)}
{\left(|x-x_0|\log\frac{1}{|x-x_0|}\right)^n}=O\left(\log\log\frac{1}{\varepsilon}\right)$$
as $\varepsilon\rightarrow 0,$ where $A(x_0, \varepsilon,
\varepsilon_0)$ is defined by (\ref{eq1**}).}
\end{proposition}

\medskip
The following statement holds.

\medskip
\begin{theorem}\label{th1_1}
{\it\, Assume that, under conditions of Theorem~\ref{th4}, instead
of the relation~(\ref{eq1E_1}) the following is true: $Q\in
FMO(x_0).$ Then there exist $\varepsilon_0>0,$
$\widetilde{\varepsilon}_0>0$ such that the relation
$$h(f(x), f(y))\leqslant \alpha\cdot
\left(\frac{\log{\frac{1}{\varepsilon_0}}}{\log{\frac{1}
{|x-x_0|}}}\right)^{\beta}\,,$$
$\beta=(A_0K_0N\omega^{\,-1}_{n-1})^{-\frac{1}{n-1}},$ holds for
some $\alpha>0,$ any $x, y\in B(x_0,
\widetilde{\varepsilon}(x_0))\cap D$ and all $f\in\frak{R}^{E_*, E,
N}_{\varphi, Q, \delta, A_0}(D)$ whenever $n\geqslant 3$
($f\in\frak{R}^{E_*, E, N}_{Q, \delta, A_0}(D)$ for $n=2$).}
\end{theorem}

\medskip
\begin{proof}
Let $\psi(t)=\frac{1}{t\,\log{\frac1t}}>0.$ Observe that, the
quantity $I(\varepsilon, \varepsilon_0)$ defined in~(\ref{eq3AB})
may be calculated as follows:
$$ I(\varepsilon,
\varepsilon_0)=\int\limits_{\varepsilon}^ {\varepsilon_0}\psi(t) dt
=\log{\frac{\log{\frac{1}
{\varepsilon}}}{\log{\frac{1}{\varepsilon_0}}}}\,.$$
In this case, by Proposition~\ref{pr2} we obtain that
$$\frac{1}{I^n(\varepsilon,
\varepsilon_0)}\int\limits_{\varepsilon<|x-x_0|<\varepsilon_0}
Q(x)\cdot\psi^n(|x-x_0|)\,dm(x)\leqslant
C_1\left(\log\log\frac{1}{\varepsilon}\right)^{1-n}\,.$$
Let $|x-x_0|\geqslant |y-x_0|.$ Now, by Lemma~\ref{lem3} with $p=1$
in (\ref{eq3.7A}), we obtain that
$$h(f(x), f(y))\leqslant \alpha\cdot
\exp\left\{-(A_0K_0N\omega^{\,-1}_{n-1})^{-\frac{1}{n-1}}\cdot
\log{\frac{\log{\frac{1}
{|x-x_0|}}}{\log{\frac{1}{\varepsilon_0}}}}\right\}$$
$$=\alpha\left(\frac{\log{\frac{1}{\varepsilon_0}}}{\log{\frac{1}
{|x-x_0|}}}\right)^{\beta}\,,$$
where $\beta=(A_0K_0N\omega^{\,-1}_{n-1})^{-\frac{1}{n-1}},$ as
required.~$\Box$

\medskip
\begin{theorem}\label{th3}
{\it\, Assume that, under conditions of Theorem~\ref{th5}, instead
of the relation~(\ref{eq1EA}) the relation $Q\in FMO(x_0)$ holds for
every $x_0\in I(P_0).$ Then
$$h(f(x), f(y))\leqslant \alpha\cdot
\left(\frac{\log{\frac{1}{\varepsilon_0}}}{\log{\frac{1}
{|x-x_0|}}}\right)^{\beta}\,,$$
$\beta=(A_0K_0N\omega^{\,-1}_{n-1})^{-\frac{1}{n-1}},$
holds for some $\alpha>0,$ $x_0=I(P_0),$ some a neighborhood $U$ of
$P_0,$ any $x, y\in U\cap D$ and all $f\in\frak{R}^{E_*, E,
N}_{\varphi, Q, \delta, A_0}(D)$ whenever $n\geqslant 3$
($f\in\frak{R}^{E_*, E, N}_{Q, \delta, A_0}(D)$ for $n=2$).}
\end{theorem}

\medskip
{\it Proof of Theorem~\ref{th2}} is based on Lemma~\ref{lem5} and is
similar to the proof of Theorem~\ref{th1_1}.~$\Box$
\end{proof}

\section{On classes of ring $Q$-mappings}

Given $\delta>0,$ $A_0>0,$ a closed sets $E_*$ in $\overline{{\Bbb
R}^n},$ $n\geqslant 2,$ a domain $D\subset {\Bbb R}^n,$ a set
$E\subset D$ (which is closed in $D$) and a Lebesgue measurable
function $Q:D\rightarrow [0, \infty]$ we denote by $\frak{F}^{E_*,
E}_{Q, \delta, A_0}(D)$ the family of all open discrete mappings
$f:D\rightarrow {\Bbb R}^n$ satisfying~(\ref{eq2*!A})--(\ref{eq8BC})
at any point $x_0\in
\partial D$ and, in addition,

\medskip
1) $C(f, \partial D)\subset E_*,$

\medskip
2) for each component $K$ of $D^{\,\prime}_f\setminus  E_*,$
$D^{\,\prime}_f:=f(D),$ there is a continuum $K_f\subset K$ such
that $h(K_f)\geqslant \delta$ and $h(f^{\,-1}(K_f), \partial
D)\geqslant \delta>0,$

\medskip
3) $f^{\,-1}(E_*)=E,$

\medskip
4) any component $G$ of $f(D)\setminus E_*$ satisfies the
relation~(\ref{eq4***}).

\medskip
Note that all the results of this work can be transferred to classes
$\frak{F}^{E_*, E}_{Q, \delta, A_0}(D).$ In particular, the
following statements are true.

\medskip
\begin{theorem}\label{th6} {\it\, Let
$D$ be a domain in ${\Bbb R}^n,$ $n\geqslant 2.$ Assume that the
conditions 1)--3) of Theorem~\ref{th4} hold.  Then there are
$\widetilde{\varepsilon}_0>0$ and $\alpha>0$ such that the relation
$$
h(f(x), f(y))\leqslant \alpha\cdot\max\{|x-x_0|^{\beta},
|y-x_0|^{\beta}\}$$
holds for any $x, y\in B(x_0, \widetilde{\varepsilon}(x_0))\cap D$
and all $f\in\frak{F}^{E_*, E}_{Q, \delta, A_0}(D),$ where
$0<\beta\leqslant 1$ is some constant depending only on $n,$ $C$ and
$A_0.$ }
\end{theorem}

\medskip
\begin{theorem}\label{th7} {\it\, Let
$D$ be a regular domain in ${\Bbb R}^n,$ $n\geqslant 2.$ Assume that
the conditions 1)--3) of Theorem~\ref{th5} hold. Then for each
$P_0\in E_D$ there exists $x_0\in \partial D$ such that
$I(P_0)=\{x_0\},$ where $I(P)$ denotes the impression of $P_0.$ In
addition, there exists a neighborhood $U$ of $P_0$ in the metric
space $(\overline{D}_P, \rho)$ such that the inequality
$$
h(f(x), f(y))\leqslant \alpha\cdot\max\{|x-x_0|^{\beta},
(|y-x_0|)^{\beta}\}$$
holds for any $x, y\in U\cap D$ and all $f\in\frak{F}^{E_*, E}_{Q,
\delta, A_0}(D),$ where $0<\beta\leqslant 1$ is some constant
depending only on $n,$ $C$ and $A_0.$}
\end{theorem}

\medskip
The proofs of Theorems~\ref{th6} and~\ref{th7} are almost word for
word the same as the proofs of Theorems~\ref{th4} and~\ref{th5}. The
only difference between these proofs is that for the Sobolev and
Orlicz-Sobolev classes we resort to establishing a modulus
inequality of type~(\ref{eq2*!A})--(\ref{eq8BC}), for which
Lemmas~\ref{thOS4.1} and~\ref{lem3A} are used. At the same time, the
classes of ring $Q$-mappings are already defined using estimates of
the distortion of the modulus, so the use of Lemmas~\ref{thOS4.1}
and~\ref{lem3A} is not required here. The rest of the reasoning in
the proofs of these theorems is word for word.

\medskip {\bf Acknowledgement}

The work was supported by the National Research Foundation of
Ukraine (Project ``Analogues of Carath\'{e}odory and Koebe-Bloch
theorems for Orlycz-Sobolev classes'', Project number 2025.02/0010).

\medskip
{\bf \noindent Victoria Desyatka} \\
Zhytomyr Ivan Franko State University,  \\
40 Velyka Berdychivs'ka Str., 10 008  Zhytomyr, UKRAINE \\
victoriazehrer@gmail.com

\medskip
\noindent{{\bf Zarina Kovba} \\
Zhytomyr Ivan Franko State University,  \\
40 Velyka Berdychivs'ka Str., 10 008  Zhytomyr, UKRAINE \\
e-mail: mazhydova@gmail.com

\medskip
\medskip
{\bf \noindent Evgeny Sevost'yanov} \\
{\bf 1.} Zhytomyr Ivan Franko State University,  \\
40 Velyka Berdychivs'ka Str., 10 008  Zhytomyr, UKRAINE \\
{\bf 2.} Institute of Applied Mathematics and Mechanics\\
of NAS of Ukraine, \\
19 Henerala Batyuka Str., 84 116 Slov'yans'k,  UKRAINE\\
esevostyanov2009@gmail.com

\end{document}